\documentclass[11pt]{article}
\usepackage{mathrsfs}
\usepackage{}

\usepackage{bbding}
\usepackage{pifont}
\usepackage{amsfonts}
\usepackage{geometry}
\usepackage{amsmath}
\usepackage{amsthm}
\usepackage{amssymb}
\usepackage{mathptmx}
\usepackage{graphicx}
\usepackage{bbding}
\usepackage{tikz}

\def\nd{\noindent}

\def\<{\leq}
\def\>{\geq}

\newtheorem{thm}{Theorem}[section]
\newtheorem{lem}{Lemma}[section]
\newtheorem{prop}{Proposition}[section]
\newtheorem{cor}{Corollary}[section]
\newtheorem{defi}{Definition}[section]
\newtheorem{rem}{Remark}[section]

\begin{document}
\title{ Birman-Murakami-Wenzl Type Algebras for Arbitrary Coxeter Systems  }
\author{Zhi Chen}
\date{}
\maketitle

\begin{abstract}
  For every Dihedral group $W_{I_{2}(k)}$ (k=2,3,4..) we define a Birman-Murakami-Wenzl type algebra $B_{I_2 (k)}(m,l)$, based on which we present a  BMW type algebra $B_{\Gamma} (m,l) $ for any Coxeter matrix $\Gamma$, as a generalization of the BMW algebras to other real types. We show the algebras $B_{I_{2}(k)}(l,m) $ have most of  algebraic conditions of the original Birman-Murakami-Wenzl algebras including semisimplicity and cellularity. Further more, we found  a $k$ dimensional representation of the Artin-Tits group $A_{I_2 (k)} $ with two parameters $l,m$ with quite neat shape in the algebra $B_{I_2 (k)}(l,m) $, which seems a natural generalization of the Lawrence-Krammer representation to Dihedral type Artin groups. We conjecture these representations are  isomorphic to the generalized Lawrence-Krammer representation define by I.Marin by using flat connections.

\end{abstract}

\section{Introduction}

For every Coxeter matrix $\Gamma$ there are a associated Coxeter group $W_{\Gamma}$, Artin group $A_{\Gamma} $ (or Artin-Tits group in some literatures ) and a Hecke algebra $H_{\Gamma}(v)$. The Birman-Murakami-Wenzl algebras (BMW algebra later on)  $B_{n} (m,l)$ can be seen as a  structure " lay above" the Hecke algebra $H_n (v)$ of $A_{n-1}$ type in many aspects.  A interesting feather of $B_n (m,l)$ is that it contains the famous Lawrence-Krammer representation (LK representation later on ) of the Braid groups $B_n$ \cite{Z} .

  When the Coxeter matrix $\Gamma$ is of simply laced type (ADE type when $\Gamma$ is of finite type), there is a very natural generalization $B_{\Gamma}(m,l ) $ of the BMW algebra and its degenerate version Brauer algebra \cite{CGW1} \cite{CFW}. The algebra $B_{\Gamma}(m,l)$ contains a sub representation of the  Artin group $A_{\Gamma}$ which can be seen a natural generalization of the LK representation \cite{CGW1} \cite{CW}.  Besides, for a Complex reflection group $G(m,1,n)$, Haring-Oldenberg defined a so called cyclotomic-BMW algebra \cite{Ha} which was studied by many authors.

  It is a  question whether there exist a BMW type algebra for any Coxeter matrix $\Gamma$. The style of generalizing in \cite{CW} only works for simply laced Coxeter matrixes. If such generalization do exist then the main problem is to deal with "multiple edges", which corresponds to dihedral groups other that $S_3 (I_2 (3))$.   If such generalization exist  then they may form a complete new system of algebra (BMW type algebra )  "laying above " Hecke algebras: also be quotient algebra of the group algebra of Artin group, they have Hecke algebra of the corresponding type as a quotient algebra and the generators of the Artin group $A_{\Gamma}$ have a degree 3 annihilating polynomial in them.

   There is a very close relationship between (generalized) LK representations and (generalized )BMW algebras for simply laced Coxeter matrixes \cite{CW} \cite{CGW1} \cite{Pa}. In \cite{Ma2}, I.Marin defined a generalized LK representation for any finite type Artin group and any braid group associated with a complex reflection group, as monodromy of certain flat connection with nice shape. It is still a question to write down those generalized LK representation explicitly.  Inspired by Marin's construction in \cite{Ch1}  we introduced a Brauer type algebra $B_{\Gamma} ( \iota)$ for any Coxeter group and any complex reflection group. A presentation of the Brauer type algebra of Dihedral type is in table 1. Those Brauer type algebras support Marin's flat connection naturally. We conjectured they can be deformed to be certain generalized BMW algebra which contain Marin's generalized LK representations.

  In a unpublished paper \cite{Ch2}, we had a uncomplete attempt  to define a deformation $B_{I_2 (k)} (m,l) $ of the Brauer type algebras $Br_{I_2 (k)}(\iota)$ of Dihedral types. Those deformation has the same dimension of the corresponding Brauer type algebra, and they contain Marin's  generalized LK representations. But these presentations have a big problem, that is, in them there are some constants determined by Marin's generalized LK representations. To know what these constants are we need to integral Marin's flat connection to get a explicit form of the generalized LK representations, which seems a quite unreasonable task.

  In this paper we find a way out, that is , we try to guess what those constants are. If those algebras are indeed natural generalizations of the BMW algebra they should have some special algebraic properties, such as having a certain type of involution which is needed for a cellular structure. (Such a involution exsit for the original BMW algebras and all Hecke algebras, and are indispensable for the construction of the Kazhdan-Lusztig polynomial)   We did some test computations about with which constants the resulted algebra would be "nice". As a result we guess those constants are actually very simple (right hand side of $(5),(6),(7),(8)$ of Definition \ref{defi:bmwodd} ; $(8),(9),(10),(11),(12)$ of Definition\ref{defi:evenBMW} ). We studied the resulted algebra and found they satisfy all our requests to be a generalization of BMW algebra. Further more, the algebra $B_{I_2 (2n+1)}(m,l)$ contains a $2n+1$ dimensional irreducible representation of the Artin-Tits group $A_{I_2 (2n+1)}$  with quite neat shape with a quite neat invariant bilinear form ( Definition \ref{defi:oddLKdefinition} ); the algebra $B_{I_2 (2n)} (\bar{\iota} )$ has two $n$ dimensional irreducible representations of the group $A_{I_2 (2n)}$ ( Theorem \ref{thm:evenLK} ). We conjecture they are isomorphic to Marin's generalized LK representations.  If  it is true then we inversed the study of these objects: we figured out the right coefficients and obtain the monodromy of Marin's flat connection without doing integration.   It isn't surprising that the odd cases $B_{I_2 (2n+1)}(m,l)$ and the even cases $B_{I_2 (2n)}(\bar{\iota}) $ are very different and need to be handled separately. The even cases are a little more complicated, certain generalized Fibonacci sequences (Definition \ref{defi:twopolynomialsequence}) appear behind those computations.

  We list the special properties of the algebra $B_{I_{2n+1}}(m,l)$ and $B_{I_{2n}}(\bar{\iota}) $ which justify them to be a reasonable generalization of BMW algebras.

  (1) They are semisimple.

  (2) They have natural cellular structures.

  (3) There are quotient algebras of the group algebra $\mathbb{Q}[m ,l^{\pm}] A_{I_{k}}$ in which the Artin group generators $\sigma_i $ have a degree 3 annihilating polynomial.

  (4) They have the corresponding Hecke algebras as a quotient.

  (5) They contain a representation of the corresponding  Artin group $A_{I_2 (k)}$ of Lawrence-Krammer type.

  (6) There is a natural involution on them.

  (7)(to be certified ) The Brauer type algebra $Br_{\Gamma}(\iota)$ defined in \cite{Ch1} support a Knizhnik-Zamolodchikov type formal connection, every representation of the Artin-Tits group $A_{\Gamma}$ could be isomorphic to monodromy of such flat connections.

Another fact giving us the confidence is that even though the computations are complicated, the algebras (Definition \ref{defi:bmwodd}, Definition \ref{defi:evenBMW} ) and the representations (Definition \ref{defi:oddLKdefinition},Theorem \ref{thm:evenLK} )  have quite simple appearances. The generalized LK representation also has a quite simple invariant bilinear form (Theorem \ref{thm:oddLKinvariantform}, \ref{thm:oddspecialprojectorexplicit} ) .  Once we find  definitions of BMW type algebras for all Dihedral Coxeter groups, immediately we can present a definition of a BMW type algebra for any Coxeter matrix $\Gamma$ (Definition \ref{defi:generaltypeBMW} ) .

\begin{tabular}{|l|l|l|}

 \multicolumn{3}{c}{Table 1 }  \\

\hline

&  $B_{G_{2k+1} }(\Upsilon) $       &   $B_{G_{2k}}(\Upsilon)$ \\

\hline

generators &  $S_0 ,S_1 , E_0 ,E_1 $   &  $S_0 ,S_1 , E_0 ,E_1 $ \\

\hline

relations  &  $1)$ $[S_0 S_1 \cdots ]_{2k+1} =[S_1 S_0 \cdots
]_{2k+1} $ ; & $1)$ $[S_0 S_1 \cdots ]_{2k} =[S_1 S_0 \cdots ]_{2k}
; $ \\

    & $2)$ $S_0 ^2 =S_1 ^2 =1 $;  & $2)$ $S_0 ^2 =S_1 ^2 =1 $;   \\
    & $3)$ $S_i E_i = E_i = E_i S_i $ for $i=0,1$;  & $3)$ $S_i E_i = E_i = E_i S_i $ for $i=0,1$; \\
     & $4)$ $E_i ^2 = \tau E_i $ for $i=0,1 $;  & $4)$ $E_i ^2 = \tau_i E_i $ for $i=0,1 $; \\
    &  $5)$ $E_0 [S_1 S_0 \cdots ]_{2i-1 } E_0$ & $5)$ $E_{0} [S_1 S_0 \cdots ]_{2i-1 } E_{0}$ \\
    &$=\mu E_0 $ for $1\leq i\leq
k$; & $= (\mu _i +\mu _{i+k}
 ) E_0 $ for $1\leq i\leq k$; \\
    & $6)$ $E_1 [S_0 S_1 \cdots ]_{2i-1 } E_1 $ & $6)$ $E_1 [S_1 S_0 \cdots ]_{2i-1} E_1 $   \\
     & $=\mu E_1 $ for $1\leq i\leq
k$;  & $= (\mu _i +\mu
_{i+k}  )E_1 $ for $1\leq i\leq k$;\\
    & $7)$ $[S_0 S_1 \cdots ]_{2k } E_0 $& $7)$ $[S_1 S_0 \cdots ]_{2k-1} E_0 $  \\
    & $= E_1 [S_0 S_1 \cdots ]_{2k } $ ; & $ = E_0 [S_1 S_0 \cdots ]_{2k-1}=
E_0$; \\
 & $8)$ $[S_1 S_0 \cdots ]_{2k } E_1 $ & $8)$ $[S_0 S_1 \cdots ]_{2k-1} E_1 $ \\
  & $= E_0 [S_1 S_0 \cdots ]_{2k }.$ &$= E_1 [S_0 S_1 \cdots
]_{2k-1}=E_1$ ;  \\
 & & $9)$ $E_1 W E_0 = E_0 W E_1 =0 $. \\

 \hline

 \end{tabular}$\\$


   Above table 1 from \cite{Ch1} are canonical presentations for dihedral Brauer algebras.  We use $G_n$ to denote the dihedral group corresponding to a regular $n$-gon , which is a rank 2 Coxeter group of type $I_2 (n)$ . The symbol $\Upsilon $ means a group of parameters. For $G_{2k+1}$, $\Upsilon$ means $\{ \mu , \tau \} $. But after some normalization we can set $\mu =1$. For $G_{2k}$ ,$\Upsilon$ means $\{ \mu _0 , \mu _1 , \tau _0 ,\tau _1    \}$. In the presentation we set $\mu _m = \mu _0 $ if $m $ is even and $\mu _m = \mu _1$ is $m$ is odd.  In relation $(9)$ for $B_{I_2 (2k)} (\Upsilon )$ , the symbol $W$ means any element in the Coxeter group $G_{2k}$.

 \section{Definition and basic structures in cases $I_2 (2n+1)$}

   Denote $\Lambda = \mathbb{Q}[m ,l^{\pm} ]$, which is a laurent polynomial ring with two variables $v,l$. On $\Lambda $ we define the following "conjugation": $f(m,l) \mapsto \bar{f}= f(-m ,l^{-1})$, for any $f(m,l)\in \Lambda$. Denote the quotient field of $\Lambda$ as $F_{\Lambda } $.

   Then $B_{I_2 (2n+1 )} (m,l)$ is the $\Lambda $ algebra with the following canonical presentation.

  \begin{defi}
  \label{defi:bmwodd}
  The dihedral BMW algebra $B_{I_2 (2n+1 )} (m ,l )$ is generated by $X_0 ,E_0 , X_1 , E_1 $ with the following relations. We set $m= v - v^{-1} $, and $\tau =\frac{l-l^{-1}}{m} +1= \frac{l- l^{-1} } {v- v^{-1}  } +1 $.
 \begin{align*}
 &1) [X_0 X_1 ... ]_{2n+1} =[X_1 X_0 ...] _{2n+1} .   &6)& E_1 [X_0 X_1 ... ]_{2i-1} E_1 = l^{-1} E_1 \  for\  1\leq i\leq n.\\
  &2) X_i -X_i ^{-1} = m (E_i -1 )\  for\  i=0,1 .  &7)&  E_0 [ X_1 ^{-1} X_0 ^{-1} ...]_{2i-1} E_0 = l E_0  \   for\  1\leq i\leq n.   \\
  &3) X_i E_i =E_i X_i = l E_i  \ for\  i=0,1.  &8) &  E_1 [X_0 ^{-1} X_1 ^{-1} ... ]_{2i-1} E_1 = l E_1 \  for\  1\leq i\leq n.  \\
 &4) E_i ^2 = \tau E_i \  for\  i=0,1 .     & 9)& [X_0 X_1 ... ]_{2n} E_0 = E_1 [X_0 X_1 ... ]_{2n} . \\
 &5) E_0 [ X_1 X_0 ...]_{2i-1} E_0 = l^{-1} E_0 \  for\  1\leq i\leq n. & 10)& E_0 [X_1 X_0 ... ]_{2n} = [X_1 X_0 ... ]_{2n} E_1 .
\end{align*}
Since the coefficients in above presentation are all in $\Lambda$, so above presentation also defines a $\Lambda -$algebra, which is denoted as $\bar{B}_{I_2 (2n+1)}(m,l) $. It is easy to see

 $B_{I_2 (2n+1)}(m,l)\cong F_{\Lambda} \otimes _{\Lambda} \bar{B}_{I_2 (2n+1)}(m,l) .$

\end{defi}

{\noindent Notations:}  By $[AB... ]_{N}$ we mean a expression whose length is $N$, and in which $A,B$ appear alternatively.
For example,  $ [ X_i X_j ...]_{4} = X_i X_j X_i X_j$, $[X_i ^{-1} X_j ^{-1}... ]_{5}= X_i ^{-1} X_j ^{-1} X_i ^{-1} X_j ^{-1} X_i ^{-1}$. Sometimes for emphasizing the last term we write $[AB...]_{2N}$ also as $[AB...B ]_{2N} $, and $[AB...]_{2N+1}$ as $[AB...A]_{2N+1}$.

  \begin{rem}
  \label{rem: Gdodd}
  (1)It is evident that above relations are lifting of the relations in table 1 for $B_{G_{2n+1}} (\Upsilon )$. The relations 1),2),3),4),9),10) are easy generalizations of relations of the original BMW algebras. The key relations are 5), 6),7) and 8). In introduction we have explained how to find these relations and the coefficient $l^{-1}, l$. If $m\neq 0$, it isn't hard to see relations 9) and 10) are consequences of 1) and 2).

  (2) There is natural quotient map $\pi : B_{I_2 (2n+1 )}(m ,l ) \rightarrow  H_{I_2 (2n+1 )}(v) $ (Definition \ref{defi:dihedralhecke}) , by sending $X_i$ to $x_i $ ,and $E_i $ to 0 for $i=0,1$. Where $H_{I_2 (2n+1 )}(v)$ is the $I_2 (2n+1)$ type Hecke algebra, and $x_0 , x_1$ is its canonical generators.

  (3) From above relations, it isn't hard to see there is an natural degree 2 anti-automorphism $\phi$ defined by sending $W$ to $W^{re}$ ,where $W$ is any word made from $X_i ^{\pm } , E_i $, and $W^{re }$ is the word reversing $W$. For example, if $W=abcd$, then $W^{re} = dcba$.

  (4) If $l=1$ ,then the algebra degenerates to the Brauer type algebra $B_{I_2 (2n+1 )}(\tau )$. In this case $m=\frac{l-l^{-1}}{\tau -1} =0 $, so relation 2) degenerate to $X_i - X_i ^{-1} =0$.

   (5) From the presentation ,it is easy to see there is a automorphism $\psi : B_{I_2 (2n+1 )}(m ,l ) \rightarrow  B_{I_2 (2n+1 )}(m,l)  $ extending the map $\psi (X_i )= X_i ^{-1} $, $\psi (E_i ) = E_i  $, $ \psi (l )= l^{-1}$, $\psi (m) = -m $.  This automorphism maybe called conjugation.

  \end{rem}

  \begin{defi} \label{defi:dihedralhecke}The Hecke algebra $H_{I_2 (2n+1)}(v)$ is a algebra with generators $x_0 , x_1$, with the following relations.

  (1) $ [x_0 x_1 ... ]_{2n+1}=[x_1 x_0 ...]_{2n+1}$,
   (2)  $( x_i +v  )(x_i - v^{-1})=0 $ for $i=0,1$.

   The Hecke algebra $H_{I_2 (2n)}(v_0 ,v_1)$ is the algebra with generators $x_0 , x_1$ with the following relations.

   (1) $[x_0 x_1 ...]_{2n} =[x_1 x_0 ... ]_{2n} $, (2) $(x_i + v_i  )(x_i - v_i ^{-1} )=0 $ for $i=0,1$.

  \end{defi}

  First we have the following lemma, which generalize the corresponding topological relations of the original BMW algebras.

  \begin{lem}
  \label{lem:bmwinterestingrelations}
  In the algebra $B_{I_2 (2n+1 )} (m ,l )$ we have:

  (1) $ E_0 E_1 = [X_1 X_0 ... ]_{2n} E_1 = [X_1 ^{-1} X_0 ^{-1} ...]_{2n} E_1 = E_0 [X_1 X_0 ...]_{2n} = E_0 [X_1 ^{-1} X_0 ^{-1} ... ]_{2n}$.

  (2) $ E_1 E_0 = [X_0 X_1 ... ]_{2n} E_0 = [X_0 ^{-1} X_1 ^{-1} ...]_{2n} E_0 = E_1 [X_0 X_1 ...]_{2n} = E_1 [X_0 ^{-1} X_1 ^{-1} ... ]_{2n}$.

  (3) $E_0 E_1 E_0 = E_0 $, $E_1 E_0 E_1 = E_1 $.

  \end{lem}

  \begin{pf} For (1), by 10) in Definition \ref{defi:bmwodd}, we have $E_0 = [X_1 X_0 ...]_{2n} E_1 [X_0 ^{-1} X_1 ^{-1} ... ]_{2n}$. Then we have
  \begin{align*}
   E_0 E_1 &= [X_1 X_0 ...]_{2n} E_1 [X_0 ^{-1} X_1 ^{-1} ... ]_{2n-1} X_1 ^{-1} E_1 = l^{-1} [X_1 X_0 ...]_{2n} E_1 [X_0 ^{-1}  X_1 ^{-1} ...]_{2n-1} E_1 \\
         &= l^{-1} l [X_1 X_0 ...]_{2n} E_1 = [X_1 X_0 ...]_{2n} E_1 .
 \end{align*}
Where the second equality sign is because of 3) in Definition \ref{defi:bmwodd}, the third one is because of 8). We also have  $E_0 = [X_1 ^{-1} X_0 ^{-1} ...]_{2n} E_1 [ X_0 X_1 ... ]_{2n} $ by 9) of Definition \ref{defi:bmwodd}. Using this euqality ,in similar way we prove  $E_0 E_1 = [X_1 ^{-1} X_0 ^{-1} ... ]_{2n} E_1 $.  The other equality sign of (1) can be proved similarly.

The proof of (2) is similar with the proof of (1).

  For (3), by using 2) of Definition \ref{defi:bmwodd} , we have

  $E_0 E_1 E_0  =  E_0 [  \frac{1}{m} (X_1 - X_1 ^{-1}) +1     ] E_0   =   [ \frac{1}{m} ( l^{-1} -l   )  + \tau  ] E_0 = E_0 $.

   Where the first equality sign is by 2) of Definition \ref{defi:bmwodd}, the second one is by 4) ,5) and 7).

 \end{pf}

 Next we compute the dimension and determine irreducible representations of $B_{I_2 (2n+1 )}(m ,l)$. The method is to first obtain a upper bound for $\dim B_{I_2 (2n+1 )}(m,l)$, then obtain a lower bound by constructing all its irreducible representations.

 \begin{lem}
 \label{lem:bmwprojectorproperty}
 For any element $x\in B_{I_2 (2n+1 )} (m,l)$, there are polynomials $l_0 (x) , l_1 (x) \in  \Lambda $ such that $ E_0 x E_0 = l_0 (x) E_0 $ and $E_1 x E_1 = l_1 (x) E_1 $.
 \end{lem}

\begin{pf}
 For convenience, let $x,y \in B_{I_2 (2n+1 )} (m ,l )$, if there is some $\alpha \in \Lambda $ such that $x= \alpha y$ , then we denote $x\sim y $.  So what we need to prove is  $E_0 x E_0 \sim E_0 $ and
 $ E_1 x E_1 \sim E_1 $.

We can suppose $x$ is a word made from $X_i ^{\pm} , E_i $. Because of relations 2), 3) and 4) of Definition \ref{defi:bmwodd}, we can suppose $x$ is a word made from $X_i , E_i $ ,$i=0,1$, and in this word $0,1$ appear alternatively in the lower indices. For example: $ x= X_0 E_1 E_0 X_1 E_0 $. We call such a word as a "alternating word ".   Essentially we only need to consider the cases that the first low indice and the last low indice are both $1$.

We do induction on the length $l(x)$.Because of 3),4) of Definition \ref{defi:bmwodd}, we only need to consider the cases when $l(x)$ is odd. Suppose we have proved $E_0 x E_0 \sim E_0 $ and $E_1 x E_1 \sim E_1 $ for any alternating word $x$ such that $l(x)\leq 2N+1$.

 Now suppose  $x$ is a word that $l(x)= 2N+3 $, we start to show $E_0  x E_0 \sim E_0 $ first.  If there is another $E_0 $ in the word $x$, then by induction we have  $E_0 x E_0 \sim E_0 $. If there are more than two $E_1$'s in $x$, also by induction $E_0 x E_0 = E_0 ...E_1 y E_1 ...E_0 \sim E_0 $ (because by induction we have $E_1 y E_1 \sim E_1 $ ). If the low index of the word $x$ start with "$0$", then by relations 3) ,4) of Definition \ref{defi:bmwodd} and induction we have $E_0 x E_0 \sim E_0$. So we are left with the cases when there are no $E_i $'s in $x$, or there is one $E_1$ in $x$. For the first case $x= [X_1 X_0 ...  ]_{2N+3} $, for the second case $x=[ X_1 X_0...]_{2l_1} E_1 [X_0 X_1 ...]_{2l_2 }$.

Cases 1.  $x= [X_1 X_0 ...]_{2l_1} E_1 [X_0 X_1 ...]_{2l_2 }$.

If one of $l_1 ,l_2$, for example $l_1 \geq n$, then $E_0 x E_0 =[X_1 X_0 ...]_{2n} E_1 [X_1 X_0 ...  ]_{2l_1 -2n} E_1 [X_0 X_1 ...]_{2l_2 } E_0 \sim E_0 $ by induction (once again since $E_1 [X_1 X_0 ...  ]_{2l_1 -2n} E_1 \sim E_1 $). So we suppose $l_1 ,l_2 <n $ . We have
\begin{align*}
&E_0 [X_1 X_0 ... ]_{2l_1 } E_1 [X_0 X_1 ... ]_{2l_2 } E_0 = E_0 [X_1 X_0 ... ]_{2n-2l_1} ^{-1} [X_1 X_0 ... ]_{2n} E_1 [X_0 X_1 ...]_{2l_2 } E_0 \\
&=  E_0 [X_1 X_0 ... ]_{2n-2l_1} ^{-1} E_0 [X_1 X_0 ... ]_{2n}  [X_0 X_1 ...]_{2l_2 } E_0
    = E_0 [X_1 X_0 ... ]_{2n-2l_1} ^{-1} E_0 [X_1 ^{-1} X_0 ^{-1} ... ]_{2n}  [X_0 X_1 ...]_{2l_2 } E_0 \\
& =  E_0 [X_1 X_0 ... ]_{2n-2l_1} ^{-1} E_0 [X_1 ^{-1} X_0 ^{-1} ... ]_{2n-2l_2 } E_0
    \sim E_0
\end{align*}
Case 2. $x = [X_1 X_0 ...]_{2N+3}$.

  If $2N+3 \leq 2n-1$, then by 5) of Definition \ref{defi:bmwodd}, we have $E_0 x E_0 \sim E_0 $.

  If $2N+3 >2n-1 $ and $2N+3 -2n \geq 2n $,  then
  \begin{align*}
  E_0 x E_0 &= E_0 E_1 [X_1 X_0 ... ]_{2N+3 -2n } E_0
  = E_0 E_1 [ X_1 X_0 ... ]_{2N+3 -4n } E_1 [X_0 X_1 ...]_{2n} \\
  & \sim E_0 E_1 [X_0 X_1 ... ]_{2n}
  = E_0 [X_1 ^{-1 } X_0 ^{-1} ... ]_{2n} [X_0 X_1 ... ]_{2n} = E_0 .
  \end{align*}
  If $2N+3 > 2n-1 $ and $2N+3 -2n < 2n $, then
  \begin{align*}
   E_0 x E_0 &= E_0 E_1 [X_1 X_0 ... ]_{2N+3 -2n } E_0
  = E_0 E_1 [X_0 X_1 ...]_{2n} [X_1 ^{-1} X_0 ^{-1}... ]_{4n-2N-3} E_0 \\
 &= E_0 [X_0 X_1 ... ]_{2n} E_0 [X_1 ^{-1} X_0 ^{-1} ... ]_{4n-2N-3} E_0 \sim E_0 .
  \end{align*}





Similarly we can prove for any alternating word $x$ with $l(x) = 2N+3$ we have $E_1 x E_1 \sim E_1$ so complete the proof.

\end{pf}

\begin{lem}
\label{lem:bmwoddrepresentationproperty}
 For any fixed $k: 0\leq k\leq 2n$, for any $A \in \{ X_0 , X_1 , E_0 , E_1 \} $, we have

$ A [...X_0 X_1 ]_{i} E_0 [X_1 X_0 ... ]_k = \sum _{j=0} ^{2n} \alpha ^{A} _{i,j} [... X_0 X_1 ]_{j} E_0 [X_1 X_0 ...]_k $.

Applying the anti-automorphism in (3) of Remark \ref{rem: Gdodd} , we have

$[... X_0 X_1  ]_k E_0 [X_1 X_0 ... ]_i A =
\sum _{j=0} ^{2n} \alpha ^{A} _{i,j} [... X_0 X_1  ]_k E_0 [ X_1 X_0 ... ]_j   $

Where $\alpha ^{A} _{i,j} \in \Lambda $ are defined in Definition \ref{defi:oddLKdefinition}.

\end{lem}

\begin{pf}

It can be proved by using Definition \ref{defi:bmwodd} and Lemma \ref{lem:bmwinterestingrelations}. Since there are no difficulty in the computations we omit the case by case check, only present proofs of two cases.
\begin{align*}
 & X_0 \cdot [...X_0 X_1 ]_{2k} E_0 [X_1 X_0 ...]_{j} = X_0 ^{2} [X_1 ...X_0 X_1 ]_{2k-1} E_0 [X_1 X_0 ... ]_{j} \\
  &= (1+ ml E_0 - m X_0  )\cdot [X_1 ...X_0 X_1  ]_{2k-1} E_0 [X_1 X_0 ... ]_{j} \\
  &= [...X_0 X_1 ]_{2k-1} E_0 [X_1 X_0 ...]_{j}
 + ml \cdot l^{-1} E_0 [X_1 X_0 ...]_{j} -m [...X_0 X_1 ]_{2k} E_0 [X_1 X_0 ...]_{j} .
\end{align*}
\begin{align*}
&E_1 \cdot [...X_0 X_1 ]_{2k} E_0 [X_1 X_0 ...]_{j}
= E_1 \cdot [X_1 ^{-1} X_0 ^{-1}... ]_{2n-2k} [X_0 ...X_1 X_0 ]_{2n} E_0 [X_1 X_0 ...]_{j}\\
&= l^{-1} E_1 \cdot [ X_0 ^{-1} X_1 ^{-1}...X_0 ^{-1} ]_{2n-2k-1} E_1 E_0 [X_1 X_0 ...]_{j}\\
&= l^{-1} \cdot l E_1 E_0 [X_1 X_0 ... ]_{j}  = [...X_0 X_1 ]_{2n} E_0 [X_1 X_0 ... ]_{j}.
\end{align*}

 We'd like to mention that if fix any $k$, $\{  [...X_0 X_1 ]_i E_0 [X_1 X_0 ...]_k  \}$ span a $2n+1$ dimensional (left )representation of the Artin group $A_{I_2 (2n+1 )}$ , whose proof have to wait for later sections. This is the origin of Definition \ref{defi:oddLKdefinition} for the generalized LK representations.
\end{pf}

\begin{cor}
\label{cor:bmwoddbasis}
 We denote

$\Sigma_1  = \{    [...X_0 X_1 ]_{l_1 } E_0 [ X_1 X_0 ... ]_{l_2 }  |  0\leq l_1 \leq 2n ; 0\leq l_2 \leq 2n  \} $

   $\Sigma_2 =  \{ [ X_1 X_0 ...]_i  |  0\leq i \leq 2n    \}   \cup \{  [X_0 X_1 ...]_j | 1\leq j\leq 2n+1     \}$,

   then the algebra $B_{I_2 (2n+1 )} (m ,l )$ is spanned by the set $\Sigma _1 \cup \Sigma _2 $. Particularly  we have:  $\dim B_{I_2 (2n+1 )} (m ,l ) \leq 2(2n+1)+(2n+1)^2 $.

\end{cor}

\begin{pf}
 Denote the linear space spanned by $\Sigma_1 \cup \Sigma_2 $ and $\Sigma_1 $ as $V$ and $V^{'}$ respectively. Direct computation shows they are closed under actions of $X_0 ,X_1  , E_0 , E_1$ from both sides. By (2) of Definition \ref{defi:bmwodd} we have $X_i ^{-1} = X_i - m (E_i -1 )$ for $i=0,1$, so $V$ and $V^{'}$ are also closed under the action of $X_i ^{-1}$'s from both sides. So we have $E_1 = [... X_0 X_1 ]_{2n} E_0 [X_1 ^{-1} X_0 ^{-1} ... ]_{2n} \in V $. Now $V$ is an ideal of the algebra, and contains all its generators so we have $V= B_{I_2 (2n+1 ) }(m ,l)$.
\end{pf}

\begin{rem}
We will show later that $\Sigma$ is a basis of $B_{I_2 (2n+1 ) } (m ,l )$. This basis is similar with the natural basis of Hecke algebras consists of positive braids lifted from irreducible words of the corresponding Coxeter group. It seems reasonable to believe the original BMW algebras and all generalized BMW algebras have such basis.

\end{rem}

 In the next section we will construct some irreducible representations for $B_{I_2 (2n+1 )} (m,l )$. By (2) of remark \ref{rem: Gdodd} we have a morphism $\pi : B_{ I_2 (2n+1 )} (m ,l ) \rightarrow H_{I_2 (2n+1 )}$ which made $H_{I_2 (2n+1 )}$ a quotient of $B_{I_2 (2n+1 )}(m ,l )$. Suppose the irreducible representations of $H_{I_2 (2n+1 )}$ are $\rho _1 , \rho _2 ,... , \rho _N $. Then we have a set of different irreducible representations of $B_{G_{2n+1}} (m ,l )$: $\rho _1 ^{'} , \rho _2 ^{'} ,... , \rho _N ^{'} $, where $\rho _i ^{'}$ is induced from $\rho _i $ by $\pi $.

\section{Generalized LK representations in cases $I_2 (2n+1 )$ }

Next we construct another irreducible representation where the actions  of $E_i$ 's aren't zero.   This is what we called the generalized LK representation of $I_2 (2n+1 )$ type. Let's  explain how do we find such a representation.  First we pretend that we have already proven the set $\Sigma $( of Corollary \ref{cor:bmwoddbasis}) is linear independent. Denote the subspace spanned by $\{ [...X_0 X_1 ]_{l} E_0 | 2n\geq l\geq 0   \}$ as $V_0 $.  By Lemma  \ref{lem:bmwoddrepresentationproperty}, we see  $V_0 $ is closed by left action of $B_{I_2 (2n+1 )} (m ,l )$, then we should have a $2n+1$ dimensional representation on $V_0 $. Then we write down the matrixes  of actions of $X_0 , X_1 , E_0 , E_1$ under the assumed basis $\{ [...X_0 X_1 ]_{l} E_0 | 2n\geq l\geq 0   \}$, which are listed in table 2.  At last through some computations we show these actions really satisfy  all relations of definition 2.1.

\begin{tabular} {|l|l|l|}

\multicolumn{3}{c}{ Table  2 } \\

\hline

 $X_0$ , &$X _0 \cdot v_0 = l v_0 $;    &         $X_1 \cdot v_0 = v_1 $;                  \\

   $X_1$ &$X _0 \cdot v_{2n} = m v_0 -m v_{2n} + v_{2n-1} $;   &    $X _1 \cdot v_{2n} = l v_{2n}$;    \\

   &$X _0 \cdot v_{2k} = m v_0 - m v_{2k} + v_{2k-1} $, &     $X _1 \cdot v_{2k} = v_{2k+1}$,     \\

   &for $1 \leq k\leq n-1 $; &    for $1\leq k\leq n-1$;    \\

    &$X _0 \cdot v_{2k+1} = v_{2k+2}  $,  &     $X _1 \cdot v_{2k+1} = v_{2k} - m v_{2k+1} + ml v_{2n}$,\\

    &for $0\leq k \leq n-1$; &  for $0\leq k \leq n-1 $.       \\

\hline

$ X_0 ^{-1} ,$ &$X_0 ^{-1} \cdot v_0 = l^{-1} v_0  $;  &   $X_1 ^{-1} \cdot v_0 = m v_0 + v_1 -m v_{2n} $; \\

  $X_1 ^{-1}$ &$X_0 ^{-1} \cdot v_{2n} = v_{2n-1} $;    &   $X_1 ^{-1} \cdot v_{2n} = l^{-1} v_{2n}$;  \\

  &$X_0 ^{-1} \cdot v_{2k} = v_{2k-1}$;  &     $X_1 ^{-1} \cdot v_{2k} = m v_{2k} + v_{2k+1} -m v_{2n}$;\\

  &$X_0 ^{-1} \cdot v_{2k+1} = -ml^{-1} v_0 + m v_{2k+1} + v_{2k+2}$;   &  $X_1 ^{-1} \cdot v_{2k+1} = v_{2k} $. \\

\hline

 $E_0$ , &$E_0 \cdot v_0 = \tau v_0 $;   &    $E_1 \cdot v_{0} = v_{2n} $;  \\

  $E_1$ &$E_0 \cdot v_{2n} = v_0 $ ;    &   $E_1 \cdot v_{2n} = \tau v_{2n} $; \\

  &$E_0 \cdot v_{2k} = v_0 $, for $1\leq k\leq n-1 $;   & $E_1 \cdot v_{2k} = v_{2n} $, for $\leq k\leq n-1$; \\

  &$E_0 \cdot v_{2k+1} = l^{-1} v_0 $.    &   $E_1 \cdot v_{2k+1} =l v_{2n} $.\\

\hline
\end{tabular}
\\

\begin{defi}
\label{defi:oddLKdefinition}
The generalized Lawrence-Krammer representation $( V_{I_2 (2n+1 )} , \rho ^{LK}  )$ of the dihedral Artin group $A_{I_2 (2n+1 )}$ is defined as follows. Suppose $X _0 , X_1$ are the canonical generators of $A_{I_2 (2n+1 )}$. First $V_{I_2 (2n+1 )} $ is a $2n+1$ dimensional linear space with a chosen basis $\{ v_0 , v_1 ,..., v_{2n} \}  $. The action of $X _0 ,X _1$ are as in the above Table 2.  For $A\in \{ X_0 , X_0 ^{-1} , X_1 , X_1 ^{-1} , E_0 , E_1  \}$, we set $\alpha ^{A} _{i,j } \in \Lambda $ for $ 0\leq i,j\leq 2n $ by the following identity.
$$ A \cdot v_i = \sum_{j=0} ^{2n} \alpha ^{A} _{i,j} v_j . $$

Matrix form of $X_0$(left one) $ X_1$ (right one )are  as follows.$\\$

$\left(
  \begin{array}{cccccccc}
    l & 0 & 0 & 0& 0 & \cdots & 0 &0\\
    0 & 0 & 1 & 0 & 0 & \cdots & 0 &0 \\
    m & 1 & -m & 0 & 0 & \cdots & 0 &0\\
    0 & 0 & 0 & 0 & 1 & \cdots & 0 &0  \\
    m & 0 & 0 & 1 & -m & \cdots & 0& 0\\
      \vdots & \vdots & \vdots & \vdots & \vdots & \ddots & \vdots &\vdots\\
      0&0&0&0&0&\cdots &0&1  \\
    m & 0 & 0 & 0 & 0 & \cdots  &1  & -m \\
  \end{array}
\right),
\left(
  \begin{array}{cccccccc}
    0& 1 & 0 & 0 & \cdots & 0 & 0 & 0 \\
    1 & -m & 0 & 0 & \cdots & 0 & 0 & ml \\
    0 & 0 & 0 & 1 & \cdots & 0 & 0 & 0 \\
    0 & 0 & 1 & -m & \cdots & 0 & 0 & ml \\
    \vdots & \vdots & \vdots & \vdots & \ddots & \vdots & \vdots & \vdots \\
    0 & 0 & 0 & 0 & \cdots & 0 & 1 & 0 \\
    0 & 0 & 0 & 0 & \cdots & 1 & -m & ml \\
    0 & 0 & 0 & 0 & \cdots & 0 & 0 & l \\
  \end{array}
\right)$ $\\$

\end{defi}

The following is the main theorem of this section,also the key step to build the algebra $B_{I_2 (2n+1 )}(m,l)$.
\begin{thm}
\label{thm:OddLKbasictheorem}
The action of $\sigma _i $($i=0,1 $) on $V_{I_2 (2n+1 )}$ defines a representation of the Artin group $A_{I_2 (2n+1 )}$ which factors through $B_{I_2 (2n+1 )}(m ,l )$.

\end{thm}

\begin{pf} We put this quite lengthy proof to the last section.
\end{pf}

\begin{rem} \label{rem:oddLKexplainasion} Since the relations $5),6),7),8)$ looks quite new let's explain more clearly.
From the proof we see the Artin relation $[X_0 X_1 ... ]_{2n+1}= [X_0 X_1 ... ]_{2n+1}$ is quite subtle. If we change the coefficients ($ l^{-1} ,l  $)  of the right side of 5),6),7),8)  of Definition \ref{defi:bmwodd}, we still have two matrixes $X_0 , X_1 $ in the same way as Lemma \ref{lem:bmwoddrepresentationproperty}, but by no means they satisfy the Artin relation. Thus the obtained algebra using new coefficients could "collapse".
\end{rem}

\begin{defi}
\label{defi:oddspecialprojector}
Here we define a set of operators which would be called special projectors. For any $0\leq i\leq 2n $, define a operator $p_i = [...X_0 X_1 ]_i E_0 [X_1 ^{-1} X_0 ^{-1} ... ]_{i} $ on $V_{I_2 (2n+1 )}$ according to above action of $X_i ,E_i$ on $V_{I_2 (2n+1 )}$.  By lemma \ref{lem:bmwprojectorproperty}, each $p_i$ is a projector to the line $\Lambda v_i $. We define a set of constants $ \alpha_{k,d}$ according to the following equality
$$ p_k ( v_d )= \alpha_{k,d} v_k . $$
\end{defi}

These operators will play a significant role in proving irreducibility of $GLK$ and constructing a invariant bilinear form for this representation. The following theorem computes those coefficients $\alpha_{k,d}$ explicitly.

\begin{thm}
\label{thm:oddspecialprojectorexplicit}
Explicitly, the action of special operator $p_i $ is as follows.
$$p_i (v_j ) = \left\{
\begin{array}{ll}
  \tau v_i , & i=j ;\\
 v_i , &if\quad i-j\equiv 0\mod 2\quad  but\quad  j\neq i ;\\
 l v_i , & if\quad  i>j \quad and\quad  i-j \equiv 1 \mod 2 ;\\
 l^{-1} v_i , &  if\quad  i<j \quad and\quad  i-j \equiv 1 \mod 2 .
\end{array}
\right.$$

\end{thm}
\begin{pf} We put this lengthy proof to the last section.
\end{pf}

\begin{thm}
\label{thm:oddLKirreducible}
The GLK representations are irreducible.
\end{thm}

\begin{pf} First, it is evident from definition of the representation that any base element $v_i $ is a generator of the representation. Now, suppose
an element $v=\sum_{i=0} ^{2n} \lambda _i v_i $ isn't a generator, then we must have $p_k (v ) =0$ for $0\leq k\leq 2n$, since $p_k$ is a projector to the line $\mathbb{C}v_k $. Then,
$0= p_k (v) =( \sum _{i=0} ^{2n} \lambda _i \alpha _{k,i} )v_k $ , so we have $\sum _{i=0} ^{2n} \lambda _i \alpha _{k,i} =0$ for $0\leq k \leq 2n$.

Now we observe in the matrix $(\alpha _{k,i })_{(2n+1) \times (2n+1 )}$, all diagonal elements are $\tau$, and the elements off diagonal are among $\{ l,l^{-1} , 1 \} $ by Theorem \ref{thm:oddspecialprojectorexplicit}. Since there are no dependence between $l$ and $\tau$, we see $Det ( (\alpha _{k,i })_{(2n+1) \times (2n+1 )} )\neq 0$ ( as an element in $F_{\Lambda }$ ). So if $v=\sum_{i=0} ^{2n} \lambda _i v_i $ isn't a generator, we must have $\lambda _i =0$ for all $i $'s.
\end{pf}
$\\$
Combining corollary 2.1, applying the Wedderburn-Artin theorem, we have the following structure theorem for $B_{I_2 (2n+1 )} (m ,l)$.

\begin{thm}
\label{thm:oddBMWsemisimpledimension}
 As an algebra over the field $F_{\Lambda }$ , the algebra $B_{I_2 (2n+1 )}(m,l )$ is semisimple.  The set $\Sigma$ is a basis of the algebra $B_{I_2 (2n+1 )}(m ,l )$, thus $\dim_{F_{\Lambda}} B_{I_2 (2n+1 )}(m ,l )= 2(2n+1)+ (2n+1)^2$.
\end{thm}

\begin{rem}
\label{rem:oddLKrealization}
So for any $0\leq k \leq 2n$, $\{ [... X_0 X_1 ]_{i} E_0 [X_1 X_0 ...]_k   \} _{i=0,...,2n} $  as a subset of $\Sigma$ is linear independent. Denote the sub $ F_{\Lambda}$ module spanned by this subset as $I_k $. By Lemma \ref{lem:bmwoddrepresentationproperty}, $I_k$ is a left ideal, and the left representation of $B_{I_2 (2n+1 )}(m,l)$ realized on this left ideal is isomorphic to the generalized LK representation  , where the isomorphism is given by
$v_i \mapsto [...X_0 X_1 ]_i E_0 [X_1 X_0 ...]_k  $ ,$i=0,1,...,2n$.
\end{rem}
Denote the elements of $\Sigma_1 \cup \Sigma_2 \subset B_{I_2 (2n+1 )}(m ,l )$  as $w_1 , ..., w_{2(2n+1) +(2n+1)^2 }$.
Suppose $\Sigma _2 $ is $\{ w_1 , ... , w_{2(2n+1)}  \}$,  $\Sigma _1 $ is $\{ w_{2(2n+1)+1 } ,..., w_{2(2n+1) +(2n+1)^2 } \} $, $w_1 = X_0 , w_2 = X_1 £¬  w_{2(2n+1)+1} = E_0 $£¬ and
 $w_{2(2n+1)+2} = [...X_0 X_1 ]_{2n} E_0 [X_1 X_0 ... ]_{2n}  $. By (1) of Lemma \ref{lem:bmwinterestingrelations}  and 9) of Definition \ref{defi:bmwodd}, actually we have
 $ w_{2(2n+1) +2} = [...X_0 X_1  ]_{2n} E_0 [ X_1 ^{-1} X_0 ^{-1} ... ]_{2n} = E_1 $.

  Since $\{ w_i \} _{i=1,..., 2(2n+1) +(2n+1)^2} $ is a basis, we see there exist a unique set of $\beta_{i,j} ^{k} (m,l) \in F_{\Lambda} $(the quotient field of $\Lambda$ ), such that the following equality holds
\begin{equation}
w_i w_j = \sum _{k=1} ^{2(2n+1)+(2n+1)^2 } \beta _{i,j} ^{k} (v,l) w_k .
\end{equation}

 We use these coefficients to define a abstract algebra. First we have the following proposition.
 \begin{prop}
 \label{prop:betaintegral}
 $ \beta_{i,j} ^{k} \in \Lambda $ for any $i,j,k$.
 \end{prop}
\begin{pf} We call an element $X \in B_{I_2 (2n+1 )}(m,l)$ as integral if $X= \Sigma_{i=1 } ^{2(2n+1) +(2n+1)^2 } a_i w_i $  such that $a_i \in \Lambda$ for any $i$.   By using Lemma \ref{lem:bmwoddrepresentationproperty}, we see if one of $w_i , w_j $ is in $\Sigma_1$ then $\beta_{i,j} ^{k} \in \Lambda $ for any $k$. Now suppose $w_i , w_j \in \Sigma_2$. It isn't hard to see for these cases, to prove $\beta_{i,j} ^{k} \in \Lambda $ we only need to prove the following four type of products are integral:
 $X_0 [X_1 X_0 ... ]_{i}  $ ($0\leq i\leq 2n $ ); $X_1 [X_1 X_0 ... ]_{i}  $ ( $ 0\leq i\leq 2n$ ); $X_1 [X_0 X_1 ... ]_{i}  $ ($0\leq i\leq 2n  $ ) and  $X_1 [X_0 X_1 ... ]_{2n+1} $ .

 This  can be proved case by case without any difficulty by using relations in Definition \ref{defi:bmwodd} and the just proved cases
 when one of $w_i , w_j$ being in $\Sigma_1$. And remember the fact $E_1 = [...X_0 X_1 ]_{2n} E_0 [X_1 X_0 ... ]_{2n} \in \Sigma_1 $.

\end{pf}

Denote the free $\Lambda$ module in $B_{I_2 (2n+1 )} (m,l) $ spanned $\Sigma_1 \cup \Sigma_2$ as $A_n $. By Proposition \ref{prop:betaintegral}, we see $A_n$ is a $\Lambda -$ algebra, and $B_{I_2 (2n+1 )} (m,l ) \cong F_{\Lambda} \otimes A_n $.

\begin{thm}
There is a natural isomorphism $\bar{B}_{I_2 (2n+1 )}(m,l)\rightarrow A_n $. So $\bar{B}_{I_2 (2n+1 )}( m,l)$ is a free $\Lambda$-module of rank $2(2n+1)+(2n+1)^2 $.
\end{thm}
\begin{pf}
First we map the generators of $\bar{B}_{I_2 (2n+1 )}(m,l)$, $X_0 , X_0 ^{-1}, X_1 , X_1 ^{-1}, E_0 , E_1 $ to
$ w_1 , w_1 -m w_{2(2n+1)+1}+m , w_2 , w_2 - m w_{2(2n+1)+2} +m , w_{2(2n+1)+1} , w_{2(2n+1)+2 }$ respectively. Since the target elements are all in $A_n$ ,and the target elements satisfy all relations in Definition \ref{defi:bmwodd}, so this map extends to a morphism $\mathcal{H}: \bar{B}_{I_2 (2n+1 )}(m,l) \rightarrow A_n$. Then we observe the target elements generate the algebra $A_n$ , so the morphism $ \mathcal{H}$ is surjective. At last the morphism induced by $\mathcal{H}$: $\bar{B}_{I_2 (2n+1 )}(m,l) \rightarrow F_{\Lambda}\otimes _{\Lambda} A_n = B_{I_2 (2n+1 )}(m,l)=F_{\Lambda}\otimes _{\Lambda} B_{I_2 (2n+1 )}(m,l)$ is nothing but the inclusion map, so $\mathcal{H}$ is also injective.
\end{pf}$\\$

Next we construct a Hermitian invariant bilinear form on the representation $V_{I_2 (2n+1 )}$. For which we need to view
$V_{I_2 (2n+1 )}$ as a free $\Lambda$ module.

 Recall the conjugation  $f(m,l) \mapsto \overline{f} = f(-m, l^{-1} ) $ for $f(m,l) \in \Lambda $  defined before the Definition \ref{defi:bmwodd}. A function $(. , . )$ : $V_{I_2 (2n+1 ) } \times V_{I_2 (2n+1 )} \rightarrow \Lambda $ is called Hermitian, if

(1) $(w ,v   ) = \overline{ (v,w)  } $ for any $v,w \in V_{I_2 (2n+1 )}$;

(2) $(f v, g w ) = f \overline{g} (v ,w ) $ for any $v,w \in V_{I_2 (2n+1 )}$ and any $f,g\in \Lambda $.

\begin{thm}
\label{thm:oddLKinvariantform}
On $V_{I_2 (2n+1 )}$ we set:  $$( v_i , v_j  )  =  \overline{ \alpha _{i,j} }$$,
  and extend it to a Hermitian form. Then this form is $A_{I_2 (2n+1 )}$ invariant. That it, for any $g \in A_{I_2 (2n+1 )}$ and any $v,w \in V_{I_2 (2n+1 )}$, we have
$ (g v , g w   )= (v,w) $. Further more, any $A_{I_2 (2n+1 )}$ invariant Hermitian form on $V_{I_2 (2n+1 )}$ is a constant multiple of this one.
\end{thm}

\begin{pf}
By Remark \ref{rem:oddLKrealization} and Lemma \ref{lem:bmwoddrepresentationproperty}, the generalized LK representation is isomorphic to the representation of $B_{I_2 (2n+1 )}(m,l)$ realized on the left ideal
$$I_0 =  \Lambda < E_0 , X_1 E_0 , X_0 X_1 E_0 ,..., [...X_0 X_1 ]_{2n} E_0  >$$
, where the isomorphism $J: V_{I_2 (2n+1 )} \rightarrow I_0 $  is given by $v_i \mapsto [...X_0 X_1 ]_{i} E_0 $.  Now by Lemma \ref{lem:bmwprojectorproperty}, we have a natural bilinear form on $I_0$ as follows:   For $x,y \in I_0 $, define $(x,y) \in F_{\Lambda}$ by the following identity.
$$ \phi \circ \psi (y) x = (x,y  ) E_0 . $$
Now for $\lambda ,\mu \in F_{\Lambda}$,
$\phi \circ \psi (\mu y ) (\lambda x )= \lambda \overline{\mu} \phi \circ \psi (y) x  $,  so
$( \lambda x ,\mu y   )= \overline {\mu} \lambda (x,y )  $, so this bilinear form is quasi linear. Next since
$$ \phi \circ \psi ( X_i x ) (X_i y )= \phi \circ \psi (x) X_i ^{-1} X_i y = \phi \circ \psi (x) y ,\quad ( i=0,1)$$
so $ (  X_i x , X_i y  )= (x,y )  $ for $i=0,1$. Since $X_0 ,X_1$ generate the Artin group $A_{I_2 (2n+1 )}$, so this bilinear form is $A_{I_2 (2n+1 )}$ invariant.

At last, since $\overline{ \alpha _{i,j } } = \alpha_{j,i}$,
$p_j (v_i ) = \alpha_{j,i} v_j $, $E_0 [X_1 ^{-1} X_0 ^{-1}... ]_{j} [...X_0 X_1 ]_{i} E_0 =  \alpha _{j,i} E_0  $,
and $\phi \circ \psi ([...X_0 X_1 ]_i E_0   ) = E_0 [X_1 ^{-1} X_0 ^{-1} ... ]_{i} $, we see the bilinear form on $I_0$ coincide with the bilinear form on $V_{I_2 (2n+1 )}$ defined in the theorem so the proof is finished.

\end{pf}

\section{Cellular structures of $B_{I_2 (2n+1 )}(m,l)$  }

Recall a cellular structure \cite{GL} on a algebra $A$ consists of the following data.

\begin{defi}
\label{defi:cellularalgebra}
A cellular algebra over $R$
is an associative algebra $A$, together with cell datum $(\Lambda
,M,C,*)$ where
\begin{itemize}
  \item (C1) $\Lambda $ is a partially ordered set and for each $\lambda \in
\Lambda$ ,$M(\lambda)$ is a finite set such that$ C:\cap _{\lambda
\in \Lambda} M(\lambda ) \times M(\lambda ) \rightarrow A $ is an
injective map with image an R-basis of A.
 \item (C2) If $\lambda \in \Lambda $ and $S,T\in M(\lambda )$,
  write $C(S,T)= C^{\lambda } _{S,T} \in A$. Then $*$ is an
  $R$-linear anti-involution of A such that $*(C^{\lambda } _{S,T} )  =C^{\lambda }
  _{T,S}$.
\item (C3) If $\lambda \in \Lambda $ and $S,T \in M(\lambda )$
  then for any element $a\in A$ we have

 $aC^{\lambda } _{S,T} \equiv \sum _{S^{'}\in M(\lambda )} r_a (S^{'} ,S) C^{\lambda} _{S^{'},T} (modA(<\lambda))
   $

   Where $r_a (S^{'} ,S)\in R$ is independent of $T$ and where $A(<\lambda
   )$ is the R-submodule of $A$ generated by $\{ C^{\mu } _{S^{''} ,T^{''} } | \mu < \lambda ; S^{''}, T^{''} \in M(\mu )
   \}.$.

\end{itemize}

\end{defi}

It is well know that there is a cellular structure on any finite type Hecke algebra. In $B_{I_2 (2n+1 )}(m ,l)$, denote the ideal generated by $E_0$ as $I_0 $ . Then the quotient algebra $B_{I_2 (2n+1 )} (m ,l ) / I_0 $ is isomorphic to the Hecke algebra $H_{I_2 (2n+1 )}(v )$. The cellular basis of $B_{I_2 (2n+1 )} (m ,l )$ ,roughly speaking consisting of two parts, a suitable lifting of the cellular basis of $H_{I_2 (2n+1 )}(v)$, and the set $\{    [...X_0 X_1 ]_{l_1 } E_0 [ X_1 X_0 ... ]_{l_2 }  |  0\leq l_1 \leq 2n ; 0\leq l_2 \leq 2n  \}$,which is a basis of $I_0$.

The anti-automorphism requested will be that one in (4) of remark 2.1. Denote this anti-automorphism as $x\mapsto *(x) $, for any $x\in B_{I_2 (2n+1 )}(m ,l )$. It is easy to see that $*(.)$ preserves the ideal $I_0 $, thus induce a degree two anti-automorphism on the Hecke algebra $H_{I_2 (2n+1 )}(v)$ which we also denote as $*$. By \cite{Ge} it is easy to see, this anti-automorphism is just the one for $H_{I_2 (2n+1 )}(v)$ needed in its cellular structure. Denote the quotient morphism from $B_{I_2 (2n+1 )} (m ,l )$ to $H_{I_2 (2n+1 )} (v)$ as $\pi$.

Suppose the cellular structure on $ H_{I_2 (2n+1 )} (v) $ as $(\Lambda ^{'}  , M^{'}  , C^{'} , *  )$, as in definition 4.1, Where the anti-automorphism $*$ is just the one as above. Then $B_{I_2 (2n+1 )}(m,l )$ have the following cellular structure
$(\Lambda , M , C , * )$.

\begin{itemize}
\item (1) The anti-automorphism $*$ is the one as above.

\item (2) $\Lambda = \Lambda ^{'} \cup { \lambda _{LK} }$. Preserve the partial order in $\Lambda ^{'}$, and for any $\lambda \in \Lambda ^{'}$, let $\lambda _{LK} < \lambda $.

\item (3) For $\lambda \in \Lambda ^{'} $, $M(\lambda ) = M^{'} (\lambda  )$;  and $M(\lambda _{LK} )= \{  0,1,2,...,2n \}$.

\item (4) For any $\lambda \in \Lambda ^{'}$, choose a total order $"<"$ in $M( \lambda)$. For $S, T \in M(\lambda )$ such that $S< T $,  set $C_{S,T}$ as any element in $\pi ^{-1} ( C^{'} _{S,T} )$.  For $T<S$, set $C_{S,T} = *(C_{T,S}) $. For any $S\in M(\lambda )$, choose any $D _{S} \in \pi ^{-1} (C^{'} _{S,S } )$, then set $C_{S,S} =  \frac {1 }{2} [ D _{S} + *(D_{S}) ]$.

\item (5) For $0\leq i ,j \leq 2n$ , set $C_{i,j} = [... X_0 X_1 ]_i E_0 [X_1 X_0 ... ]_j $.

\end{itemize}

\begin{thm}
\label{thm:oddBMWcellular}
 Above data $(\Lambda , M , C , * )$ defines a cellular structure on $B_{I_2 (2n+1 )} (m,l )$.
\end{thm}

\begin{pf}
$(C1), (C2)$ of definition 4.1 are easy to show. For $(C3)$, the statement for $\lambda \in \Lambda ^{'}$ follows from similar statement for cellularity of the data $(\Lambda ^{'} , M^{'} , C^{'} , *  )$. Now we consider the statement for $\lambda_{LK}$.

First by definition of $*$, we have $*( C_{i,j} ) = *(  [... X_0 X_1 ]_i E_0 [X_1 X_0 ... ]_j   )=  [... X_0 X_1 ]_j E_0 [X_1 X_0 ... ]_i = C_{j,i} $. Then lemma 2.3 shows there are polynomials  $\alpha ^{X_0 } _{i,j},  \alpha ^{X_1 } _{i,j }$  without dependence on $k$ such that:

$X_0 ( C_{i,k} ) = \sum _{j=0 } ^{2n}  \alpha ^{X_0 } _{i,j} C_{j,k} $; $X_1 ( C_{i,k} ) = \sum _{j=0 } ^{2n}  \alpha ^{X_1 } _{i,j} C_{j,k} $.

So the cellularity is proved.
\end{pf}

\section{Definition and basic structures in cases $I_2 (2n)$}

Like the case of Hecke algebras, $I_2 (2n)$ type generalized BMW algebra involve more parameters. We set $\bar{\Lambda} ^{'} = \mathbb{Q}[v_0 ^{\pm}, v_1 ^{\pm} , l_0 ^{\pm} , l_1 ^{\pm} ]$.  And let $\Lambda ^{'}$ be the ring obtained by localizing $\bar{\Lambda} ^{'} $ at $\{ m_i = v_i - v_i ^{-1} \} _{ i=0,1 }  $.  On $\Lambda ^{'}$ we also define a conjugation $\bar {(.)}$ by $\overline {f(v_0 , l_0 , v_1 ,l_1 )} = f(v_0 ^{-1} , l_0 ^{-1} , v_1 ^{-1} , l_1 ^{-1} )$. Denote the quotient field of $\Lambda^{'}$ as $F_{\Lambda^{'}}$.  The type $I_2 (2n)$ generalized BMW algebra are defined as follows.

\begin{defi}
\label{defi:evenBMW}
The type $I_2 (2n)$ generalized BMW algebra  $B_{I_2 (2n)}(\bar{\iota})$ is generated by $X_0 , E_0 , X_1 , E_1$ submitting to the following relations. Let $\Delta = [X_0 X_1 ... ]_{2n} $. Let $m_i = v_i - v_i ^{-1}$ and $\tau _i = \frac{l_i - l_i ^{-1}}{m_i} +1 $ for $i=0,1$.
\begin{align*}
& 1) [X_1 X_0 ...]_{2n}=[X_0 X_1 ...]_{2n}.  &7)& [X_1 X_0 ... ]_{2n-1}E_0 =E_0 [X_1 X_0 ...]_{2n-1}.\\
& 2) X_i - X_i ^{-1} = m_i (E_i -1 ), i=0,1. &8)& E_0 [X_1 X_0...]_{4k+1}E_0 =\frac{m_1}{m_0}(v_0 ^{-1}+l_0 ^{-1})E_0 ,4k+1<n.\\
& 3) X_i E_i = E_i X_i = l_i E_i ,i=0,1. &9)&  E_0 [X_1 X_0...]_{4k+3}E_0 =(v_0 ^{-1}+l_0 ^{-1})E_0 ,4k+3<n. \\
& 4) E_i ^2 =\tau _i E_i ,i=0,1.  &10)& E_1 [X_0 X_1 ...]_{4k+1}E_1 = \frac{m_0}{m_1}(v_1 ^{-1} + l_1 ^{-1})E_1 , 4k+1<n. \\
& 5) [X_0 X_1 ...]_{2n-1}E_1 = E_1 [X_0 X_1 ...]_{2n-1}.  &11)& E_1 [X_0 X_1 ...]_{4k+3}E_1 =(v_1 ^{-1} +l_1 ^{-1})E_1 , 4k+3<n. \\
&6) E_0 W E_1 = E_1 W E_0 = 0 , for\ any\   &12)& \Delta E_i = l_i v_i ^{-1} E_i , i=0,1. \\
&word\ W\   made\  from\  X_i ^{\pm} , E_i.
\end{align*}

Similar with Definition \ref{defi:bmwodd},since all coefficients in above relations are in $\Lambda ^{'}$, it also defines a $\Lambda ^{'}$ algebra which we denote as $\bar{B}_{I_2 (2n)} (\iota ) $. We have $B_{I_2 (2n)}(\iota) \cong F_{\Lambda ^{'}} \otimes _{\Lambda^{'}} \bar{B}_{I_2 (2n)} (\iota) .$

\end{defi}

\begin{rem}
(1) Relation $6)$ corresponds to $9)$ for $B_{G_{2k}}(\iota)$ in table 2.

(2) As in the cases of $B_{I_2 (2n+1)}(m ,l) $, the coefficients in the right sides of $8),9),10),11)$ are first conjectured based on some test computations ,requesting the resulted algebra satisfy some algebraic properties, then certified to be "really nice" through further investigations. Compare to the $I_{2}(2n+1)$ cases the algebra has two systems ($v_0 , l_0 ; v_1 , l_1 $  )  of parameters corresponds to the fact that the set of reflections of $W_{I_2 (2n)}$ fall into two conjugacy classes. For special cases such that $v_0 =v_1 =v ; l_0 = l_1 =l$, then the coefficients on the right side of $8) , 9) ,10) ,11)$ all become $v^{-1} + l^{-1} $.

(3) There is also a natural quotient map $\pi : B_{I_2 (2n)}(\bar{\iota})\rightarrow H_{I_2 (2n)}(v_1 ,v_2) $, by sending $X_i $ to $x_i$, and $E_i $ to $0 $ for $i=0,1$. Where $H_{I_2 (2n)}(v_1 ,v_2)$ is the $I_2 (2n)$ type Hecke algebra (Definition \ref{defi:dihedralhecke} ) , and $x_0 , x_1$ are its canonical generators.

(4) If $l_{i} =1 , v_{i} =1$ for $i=0,1$ ,then the algebra $B_{I_2 (2n)}(\bar{\iota})$ degenerate to the Brauer type algebra $B_{G_{2n}} (\iota) $ for suitable parameters $\iota$.

(5) It is evident from the definition that the map $X_i \mapsto X_i ; E_i \mapsto E_i $ ($ i=0,1$ ) extends to a degree 2 anti-automorphism of the algebra $B_{I_2 (2n)} (\bar{\iota})$, which is denoted by $ J$.

\end{rem}

As Lemma 2.2 for the algebra $\bar{B}_{I_2 (2n+1)}(m,l)$, we have the following lemma.

\begin{lem}
\label{lem:evenBMWprojectorproperty}
For any element $x\in \bar{B}_{I_2 (2n)}(\bar{\iota })$, there are polynomials $\alpha_0 (x) ,\alpha_1 (x) \in \Delta^{'}$ such that $ E_0 x E_0 = \alpha_0 (x) E_0 $ and $E_1 x E_1 = \alpha_1 (x) E_1 $.
\end{lem}

\begin{pf}
The proof is the same as in the proof of Lemma 2.2.
\end{pf}

\begin{cor}
\label{cor:evenBMWbasis}
The algebra $B_{I_2 (2n)}(\bar{\iota})$ is spanned by the subset $\Sigma_0 \cup \Sigma_1 \cup \Sigma^{'}$, where

$\Sigma_0 =\{ [...X_0 X_1 ]_{i_1 } E_0 [X_1 X_0 ... ]_{i_2 } |   0\leq i_1 , i_2 \leq n-1 \} $,

$\Sigma_1 = \{ [...X_1 X_0 ]_{i_1} E_1 [X_0 X_1 ... ]_{i_2} | 0\leq i_1 , i_2 \leq n-1 \}$,

$\Sigma^{'} = \{ [X_0 X_1 ...]_{i} | 0\leq i\leq 2n \} \cup \{ [X_1 X_0 ...]_{j} | 1\leq j\leq 2n-1 \} $.

Especially, we have $\dim_{F_{\Lambda^{'}}} B_{I_2 (2n)}(\bar{\iota}) \leq 4n+ 2 n^2$.
\end{cor}

\begin{pf}
Denote the $\Delta ^{'}$ submodules in $B_{I_2 (2n)}(\bar{\iota})$ spanned by $\Sigma_0 \cup \Sigma_1 \cup \Sigma^{'}$ , $\Sigma_0 $ and $\Sigma_1 $ as $V$, $V_0$ and $ V_1$ respectively. Since $V$ contains the generators $X_i ,E_i$,$i=0,1$, if we can prove $V$ is an ideal then the corollary follows. The main part of the proof parallels the proof in corollary 2.1, we only mention the place where relation $12)$ of Definition 5.1 is used.

The action of $X_0 $ on  elements $[...X_0 X_1]_{2n-1} E_0 [X_1 X_0 ...]_{j} (1\leq j\leq 2n-1) $ is as follows.
\begin{align*}
&X_0 [X_1 ... X_0 X_1 ]_{2n-1} E_0 [X_1 X_0 ... ]_{j} =
l_0 ^{-1} X_0 [X_1 ...X_0 X_1 ]_{2n-1} X_0 E_0 [X_1 X_0 ...]_{j} =
l_0 ^{-1} [...X_1 X_0 ]_{2n+1} E_0 [X_1 X_0 ...]_{j} \\
&= l_0 ^{-1} [X_1 ^{-1} X_0 ^{-1}...  X_1 ^{-1} ]_{2n-1 } \Delta E_0 [X_1 X_0 ... ]_{j}
= v_0 ^{-1} [X_1 ^{-1} ...X_0 ^{-1} X_1 ^{-1} ]_{2n-1} E_0 [X_1 X_0 ... ]_{i} .
\end{align*}
Where the first equlity sign is by relation $3)$ of Definition 5.1, and the forth equality sign is by relation $12)$.

By substituting $ X_i ^{-1}$ in the resulted term with $X_i - m_i E_i +m_i  $, it is easy to show this term is in $V$.
Similarly we can prove

 $[...X_0 X_1 ]_{j} E_0 [X_1 X_0 ... ]_{2n-1}\cdot X_0 ,  X_1 \cdot [...X_1 X_0 ]_{2n-1} E_1 [X_0 X_1 ...]_{j} , [... X_1 X_0 ]_{j} E_1 [X_0 X_1 ... ]_{2n-1} \cdot X_1  \in V $.

Above proof in fact implies that $ V_0 $ and $V_1 $ are ideals.
\end{pf}\\

In the following we construct two $n$ dimensional representations of $B_{I_2 (2n)} (\bar{\iota})$, which are to be called the generalized Lawrence-Krammer representations of the Artin group $A_{I_2 (2n)}$. The expression is more complex than the cases of $B_{I_2 (2n+1)}(m,l)$, yet much simpler than we anticipated. Since all complex computations turn out to be amalgamated in two polynomial sequences (essentially determined by one generalized Fibonacci sequence ). That is one reason made us believe that these representations are interesting and could be the correct generalization of LK representations of types $I_2 (2n)$.  The following definition introduces two polynomial sequences being important for construction of generalized LK representations and the algebra $B_{I_2 (2n)}(\bar{\iota } )$.

\begin{defi}
\label{defi:twopolynomialsequence}
The two related polynomial sequences $a_0 , a_1 , a_2 ...$ and $b_0 , b_1 , b_2 ...$ are defined inductively as follows.

(1) $a_0 =1 $, $ a_1 = m_1$ ; $b_0 =1 $, $b_1 = m_0 $;

(2) $a_{2i+1} = m_1 b_{2i} + b_{2i-1 }  $, $a_{2i} =b_{2i} $;

(3) $b_{2i} = -v_0 ^{-1} a_{2i-1}  $, $b_{2i-1} = m_0 a_{2i-1} + a_{2i-3} $.

Let $\digamma : \Delta^{'} \rightarrow \Delta^{'} $ be the automorphism by sending $l_0 , v_0 , l_1 , v_1 $ to $l_1 , v_1 , l_0 ,v_0  $ respectively. And set $a^{'} _i = \digamma (a_i  )$ and $b^{'} _i = \digamma (b_i )$ for all $i$.
\end{defi}

\begin{rem}
(1)Concretely, $(b_2 , b_3 ,b_4 ,... )= (-m_1 v_0 ^{-1} , -m_0 m_1 v_0 ^{-1} +m_1 , m_1 ^2 v_0 ^{-2} -m_0 v_0 ^{-1} ,... ) $,  $(a_2 , a_3 , a_4 ,... )= (-m_1 v_0 ^{-1} , -m_1 ^2 v_0 ^{-1} +m_0 , m_1 ^2 v_0 ^{-2} - m_0 v_0 ^{-1} ,... )$.

(2) It isn't hard to see those two sequences are determined by the subsequence $a_1 , a_3 ,..., a_{2i+1} ,... $ and we have $a_{2i+1} = m_1 b_{2i} +b_{2i-1} = -m_1 v_0 ^{-1} a_{2i-1} + (1-m_0 v_0 ^{-1} ) a_{2i-3} $. So this subsequence is a generalized Fibonacci sequence.
\end{rem}

For convenience, we set $\lambda_{2i+1} = \frac{m_1}{m_0} (v_0 ^{-1} + l_0 ^{-1} )  $ for $i\equiv 0\mod 2 $, and $\lambda_{2i+1} = v_0 ^{-1} + l_0 ^{-1} $ for $i\equiv 1\mod 2$.

\begin{lem}
\label{lem:polynomialsequencebasisequality}
(1) For $k\in \mathbb{Z} ^{>0}$, $(l_0 ^{-1} + v_0 ^{-1} )  a_{2k+1} = \sum_{i=0} ^{k} m_0 a_{2k-2i} \lambda_{2i+1} . $

(2) $\sum_{i=1} ^{2k} a_{4k+1-2i} l_0 \lambda_{2i-1} + \frac{l_0 +v_0 }{m_0 } a_{4k+1}=\frac{ m_1}{m_0 } (l_0 +v_0 )   $.

(3) $\sum_{i=1} ^{2k+1} a_{4k+3-2i} l_0 \lambda_{2i-1} + \frac{ l_0 + v_0  }{m_0 } a_{4k+3} =  (l_0 + v_0 ) $.
\end{lem}

\begin{pf}
For (1), by relations in Definition 5.2, we have
\begin{align*}
a_{2k+1} &= m_1 b_{2k} + b_{2k-1} = m_1 a_{2k} + m_0 a_{2k-2} + a_{2k-3}
=m_1 a_{2k} + m_0 a_{2k-2} + m_1 b_{2k-4} +b_{2k-5} \\
&=...=m_1 a_{2k} + m_0 a_{2k-2} + m_1 a_{2k-4} + m_0 a_{2k-6} +... .
\end{align*}

From which the statement follows. For (2), we first prove the following identity $(2)^{'}$ by induction.
$$(2)^{'}. a_1 m_0 +a_3 m_1 +...+a_{4k-1} m_1 =m_1 v_0 +(-v_0 ) a_{4k+1} .$$

When $k=1$, the left hand side is
\begin{align*}
a_1 m_0 + a_3 m_1 &=m_1 m_0 + (-m_1 ^{2} v_0 ^{-1} +m_0  )m_1 =
m_1 v_0 + (-m_1 v_0 + 2m_0 m_1 - m_1 ^{3} v_0 ^{-1}    ) \\
&= m_1 v_0 + (-v_0 ) a_{5} .
\end{align*}
Suppose $(2)^{'}$ have
been proven for $ k\leq K$, now for $k=K+1$, the left hand side is
\begin{align*}
&a_1 m_0 + a_3 m_1 + ...+ a_{4K-1} m_1 + a_{4K+1} m_0 + a_{4K+3} m_1 \\
&= m_1 v_0 +(-v_0 ) a_{4K+1}  + a_{4K+1} m_0 + a_{4K+3} m_1 = m_1 v_0 + (-v_0 ) [ -m_1 v_0 ^{-1} a_{4K+3} +( 1-m_0 v_0 ^{-1} ) a_{4K+1} ]\\
&=m_1 v_0 + (-v_0 ) a_{4K+5} .
\end{align*}

So by $(2)^{'}$, the left hand side of $(2)$ is
\begin{align*}
 &\frac{l_0  }{m_0 } ( v_0 ^{-1} + l_0 ^{-1}  ) [ a_1 m_0 + a_3 m_1 +...+ a_{4k-1} m_1 ] + \frac{l_0 + v_0  }{m_0 } a_{4k+1} \\
 &= \frac{l_0  }{m_0 } ( v_0 ^{-1} + l_0 ^{-1}  ) [ m_1 v_0 + (-v_0 ) a_{4k+1} ] + \frac{l_0 + v_0  }{m_0 } a_{4k+1} = \frac{m_1 }{m_0} (l_0 + v_0 ).
\end{align*}
Similarly for $(3)$ we first prove the following
$$(3)^{'}.  a_1 m_1 + a_3 m_0 +...+ a_{4k+1} m_1 = m_0 v_0 - v_0 a_{4k+3} .$$
When $k=1$, the left hand side is
\begin{align*}
&a_1 m_1 a_3 m_0 + a_5 m_1 = a_1 m_1 + a_3 v_0 + (-v_0 ) [-m_1 v_0 ^{-1} a_5 + ( 1-m_0 v_0 ^{-1} ) a_3   ] \\
 &= m_1 m_1 + (-m_1 ^{2} v_0 ^{-1} +m_0  ) v_0 + (-v_0 ) a_7 = m_0 v_0 - a_7 v_0  .
\end{align*}
Suppose we have proven $(3)^{'}$ for $k\leq K$, for $k=K+1$,  the left hand is
\begin{align*}
 &a_1 m_1 + a_3 m_0 +...+ a_{4K+5} m_1 = a_1 m_1 + a_3 m_0 + ...+ a_{4K-1} m_0 + a_{4K+1} m_1 + a_{4K+3} m_0 + a_{4K+5} m_1   \\
&=  m_0 v_0 -v_0 a_{4K+3 } + a_{4K+3} m_0 + a_{4K+5} m_1 = m_0 v_0 - v_0 [  -m_0 v_0 ^{-1} a_{4K+5} + (1-m_0 v_0 ^{-1} ) a_{4K+3}   ]   \\
&=     m_0 v_0 - v_0 a_{4K+7}.
\end{align*}
So the induction is completed. Now using $(3)^{'}$ ,the left hand side of $(3)$ is
\begin{align*}
 &\frac{l_0 }{m_0 } (l_0 ^{-1} + v_0 ^{-1}  ) ( a_1 m_1 + a_3 m_0 +...+ a_{4k+1} m_1   ) + \frac{l_0 + v_0  }{m_0 }  a_{4k+3} \\
&=\frac{l_0 }{m_0 } (l_0 ^{-1} + v_0 ^{-1}  ) ( m_0 v_0 - v_0 a_{4k+3}  ) + \frac{l_0 + v_0  }{m_0 }  a_{4k+3} =  l_0 + v_0 .
\end{align*}

 So the proof is finished.

\end{pf}

\begin{thm}
\label{thm:negativegeneratoraction}

(1) $[... X_0 ^{-1} X_1 ^{-1} ]_{2k} E_0  = \sum _{i=0 } ^{2k }   b_{2k-i } [... X_0 X_1  ]_{i} E_0    $;

(2) $ [... X_0 ^{-1} X_1 ^{-1} ]_{2k+1} E_0 = \sum _{i=0} ^{2k+1}  a_{2k+1-i} [...X_0 X_1 ]_{i} E_0  $.

\end{thm}

\begin{pf}
We prove the statements by inductions.  First if $k=0$, then

The left hand side $= X_1 ^{-1} E_0 = (X_1 - m_1 E_1 + m_1 ) E_0 = m_1 E_0 + X_1 E_0 =a_1 E_0 + a_0 E_1 = $ the right hand side. When $k=0$, the statement $(2)$ is evident.   Suppose we have proven when $k\leq K$ statements $(1)$  and $(2)$ are true. Then we consider the cases when $k= K+1$.  First for $(1)$,
$[... X_0 ^{-1} X_1 ^{-1} ]_{2K+2} E_0 = (X_0 -m_0 E_0 +m_0 ) \cdot [...X_0 ^{-1} X_1 ^{-1} ]_{2K+1} E_0 $ is the sum of three terms as follows.
\begin{align*}
A:\  &X_0 \cdot (   \sum _{i=0} ^{2K+1}  a_{2K+1-i} [...X_0 X_1 ]_{i} E_0  )\\
&=X_0 \cdot (\sum_{i=0} ^{K} a_{2K+1-2i} [...X_0 X_1 ]_{2i} E_0  ) +
X_0 \cdot ( \sum_{i=0}^{K} a_{2K-2i} [...X_0 X_1 ]_{2i+1} E_0  )\\
&=a_{2K+1} l_0 E_0 + \sum_{i=1} ^{K} a_{2K+1-2i} (X_0 )^{2} \cdot [...X_0 X_1 ]_{2i-1} E_0 +
\sum_{i=0} ^{K} a_{2K-2i} [...X_0 X_1]_{2i+2} E_0 \\
&=a_{2K+1} l_0 E_0 + \sum_{i=1} ^{K} a_{2K+1-2i} ( 1+m_0 l_0 E_0 -m_0 X_0 ) \cdot [...X_0 X_1 ]_{2i-1} E_0 +
\sum_{i=0} ^{K} a_{2K-2i} [...X_0 X_1]_{2i+2} E_0 \\
&=a_{2K+1} l_0 E_0 + \sum_{i=1} ^{K} a_{2K+1-2i} [...X_0 X_1 ]_{2i-1} E_0 +
\sum_{i=1} ^{K} m_0 l_0 a_{2K+1-2i} \lambda_{2i-1} E_0 \\
&-\sum_{i=1} ^{K} m_0 a_{2K+1-2i} [...X_0 X_1 ]_{2i} E_0
+\sum_{i=0} ^{K} a_{2K-2i} [...X_0 X_1 ]_{2i+2} E_0   .
\end{align*}
\begin{align*}
B:\ &-m_0 E_0 (  \sum _{i=0} ^{2K+1}  a_{2K+1-i} [...X_0 X_1 ]_{i} E_0   ) \\
& =
-m_0 a_{2K+1} \tau_0 E_0 - \sum_{i=1} ^{K} a_{2K+1-2i} m_0 E_0 [...X_0 X_1 ]_{2i} E_0
-\sum_{i=0} ^{K} a_{2K-2i} m_0 E_0 [...X_0 X_1 ]_{2i+1} E_0 \\
&=-m_0 a_{2K+1} \tau_0 E_0 - \sum_{i=1} ^{K} a_{2K+1-2i} m_0 l_0 \lambda_{2i-1} E_0
-\sum_{i=0} ^{K} a_{2K-2i} m_0 \lambda_{2i+1} E_0 .
\end{align*}

$C:$ $m_0 [...X_0 ^{-1} X_1 ^{-1}]_{2K+1} E_0 = \sum_{i=0 } ^{2K+1} m_0 a_{2K+1-i} [...X_0 X_1 ]_{i} E_0  $.

The coefficient of $E_0$ is

$l_0 a_{2K+1} +\sum_{i=1} ^{K} a_{2K+1-2i} m_0 l_0 \lambda_{2i-1} -m_0 a_{2K+1} \tau_0
-\sum_{i=1} ^{K} a_{2K+1-2i} m_0 l_0 \lambda_{2i-1} $

$-\sum_{i=0} ^{K} a_{2K-2i} m_0 \lambda_{2i+1} +m_0 a_{2K+1}
= -v_0 ^{-1} a_{2K+1} = b_{2K+2}   $;

the coefficient of $[...X_0 X_1 ]_{2i-1} E_0 $ is $a_{2K+1-2i} +m_0 a_{2K-2i+2} =b_{ ( 2K+2  ) - ( 2i-1 ) } $;

the coefficient of $[...X_0 X_1 ]_{2i} E_0 $ for $K\geq i>0$ is

$a_{2K-2i+2} -m_0 a_{2K+1 -2i } +m_0 a_{2K+1 -2i} =a_{2K-2i+2}= b_{2K+2-2i} $;

at last the coefficient of $[...X_0 X_1 ]_{2K+2} E_0 $ is $a_0 = 1= b_0 $. So the induction for statement $(1)$ is proved. For statement $(2)$, by induction (we also need (1) for $K+1$ just proved) we have
\begin{align*}
&X_1 ^{-1} \cdot [...X_0 ^{-1} X_1 ^{-1} ]_{2K+2}E_0 =
( X_1 -m_1 E_1 +m_1 )\cdot ( \sum_{i=0} ^{2K+2} b_{2K+2-i} [...X_0 X_1  ]_{i} E_0  )\\
&=\sum_{i=0} ^{K+1 } b_{2K+2-2i} X_1 \cdot [...X_0 X_1 ]_{2i} E_0
+ \sum_{i=0} ^{K} b_{2K+2 -2i-1}X_1 \cdot [...X_0 X_1 ]_{2i+1} E_0 +\sum_{i=0} ^{2K+2} m_1 b_{2K+2-i} [...X_0 X_1 ]_{i}E_0\\
&= \sum_{i=0} ^{K+1 } b_{2K+2-2i}  [...X_0 X_1 ]_{2i+1} E_0
+ \sum_{i=0} ^{K} b_{2K+2 -2i-1} (1+m_1 l_1 E_1 -m_1 X_1 )\cdot [...X_0 X_1 ]_{2i} E_0  \\
&+ \sum_{i=0} ^{2K+2} m_1 b_{2K+2-i} [...X_0 X_1 ]_{i}E_0 .
\end{align*}

So the coefficient of $[...X_0 X_1 ]_{2i}E_0 $ is $b_{2K+1-2i} +m_1 b_{2K+2 -2i} = a_{2K+3-2i}$;

the coefficient of $[...X_0 X_1 ]_{2i+1} E_0 $ is $b_{2k+2-2i } -m_1 b_{2K+1-2i} +m_1 b_{2K+2-2i-1} = b_{2K+2-2i}
=a_{2K+2-2i}$.  So the proof is finished.

\end{pf}

\begin{thm}
\label{thm:evenBMWinvolution}
On the algebra $B_{I_2 (2n)} (\bar{\iota })$ ,the correspondence $l_i \mapsto l_i ^{-1} £»v_i \mapsto v_i ^{-1} $;
$X_i \mapsto X_i ^{-1} ; E_i \mapsto E_i $ $(i=0,1)$ extends to a degree 2 isomorphism. We denote this automorphism as $\Psi $.

\end{thm}

\begin{pf}
It is evident above correspondence keeps relations in Definition 5.1 except $(8)~(11)$. To prove relations $(8), (9)$ are also kept we only need to show that

$(8)^{'}$ $E_0 [X_1 ^{-1} X_0 ^{-1} ... ]_{4k+1} E_0 = \frac{m_1 }{m_0 } (v_0 + l_0 ) E_0  $;
$(9)^{'}$ $E_0 [X_1 ^{-1} X_0 ^{-1} ... ]_{4k+3} E_0 = (v_0 + l_0 ) E_0 $.

Now we have  \begin{align*}
 &E_0 [X_1 ^{-1} X_0 ^{-1} ... ]_{4k+1} E_0 = E_0 (\sum_{i=0} ^{4k+1} a_{4k+1-i} [...X_0 X_1 ]_{i}   )E_0 \\
&= a_{4k+1} \tau_0 E_0 + \sum_{i=1} ^{2k} a_{4k+1-2i} E_0 [...X_0 X_1 ]_{2i} E_0
+\sum_{i=0} ^{2k} a_{4k+1-2i-1} E_0 [ ...X_0 X_1 ]_{2i+1} E_0 \\
&=[ a_{4k+1} (\frac{l_0 -l_0 ^{-1} }{m_0 } +1 ) + \sum_{i=1} ^{2k} a_{4k+1-2i} l_0 \lambda_{2i-1} + \frac{l_0 ^{-1} +v_0 ^{-1}}{m_0 } a_{4k+1}  ] E_0 \\
&= \frac{m_1 }{m_0 } (v_0 + l_0  )E_0  .
\end{align*}
  Where the first equality sign is by Theorem 5.1, the last equality sign is by (2) of Lemma 5.2.  By using (3) of Lemma 5.2, we can prove $(9)^{'}$ in similar way.
And we can prove relations $(10) ,(11)$ are kept similarly also.
\end{pf}

\section{Generalized  LK representations in cases $I_2 (2n)$   }

 In the following table we define a $n$ dimensional representation of $B_{I_2 (2n)}(\iota )$.

\begin{tabular} {|l|l|l|}
\multicolumn {3}{c}{Table 3} \\
\hline
$X_i \cdot u_j $  &  $X_0 \cdot u_0 = l_0 u_0$ ;  & $ X_1 \cdot u_0 = u_1   $; \\

 $(j<n-1)$  &  $X_0 \cdot u_{2i-1} = u_{2i}$ ;  &  $X_1 \cdot u_{2i-1} =u_{2i-2} -m_1 u_{2i-1}  $;   \\
           &  $X_0 \cdot u_{2i} = m_0 l_0 \lambda_{2i-1} u_0 + u_{2i-1}-m_0 u_{2i} $   & $X_1 \cdot u_{2i} =u_{2i+1}$ ; \\

\hline

$X_i \cdot u_{n-1}$ & $ X_0 \cdot u_{n-1} = \sum_{i=0} ^{n-1} v_0 ^{-1}  a_{n-1-i}u_i   $   & $X_1 \cdot u_{n-1} = u_{n-2} -m_1 u_{n-1} $  \\

   &  for $n\equiv 0 \mod 2\mathbb{Z}$;   & for $n\equiv 0 \mod 2\mathbb{Z} $;   \\

    &  $X_0 \cdot u_{n-1} = m_0 l_0 \lambda_{n-2} u_0 + u_{n-2} -m_0 u_{n-1}   $  &  $X_1 \cdot u_{n-1} =
    \sum_{i=0} ^{n-1} v_0 ^{-1} b_{n-1-i} u_i    $;   \\

    & for $n\equiv 1 \mod 2\mathbb{Z}$; & for $n\equiv 1\mod 2\mathbb{Z}$ ;     \\

\hline

$X_i ^{-1} \cdot u_j $ &   $X_0 ^{-1} u_0 = l_0 ^{-1} u_0  $;   & $X_1 ^{-1} u_0 = m_1 u_0 + u_1 $; \\
$(j<n-1 )$  &  $X_0 ^{-1} \cdot u_{2i-1} = -m_0 \lambda_{2i-1} u_0 + m_0 u_{2i-1} +u_{2i}  $;    &
 $X_1 ^{-1} \cdot u_{2i-1} = u_{2i-2} $ ; \\
       & $X_0 ^{-1} \cdot u_{2i} =u_{2i-1} $;      &  $X_1 ^{-1}\cdot u_{2i} = m_1 u_{2i} + u_{2i+1} $;  \\

\hline
$X_i ^{-1}  \cdot u_{n-1} $ &  $X_0 ^{-1} \cdot u_{n-1} =(v_0 ^{-1} a_{n-1} -m_0 \lambda_{n-1} ) u_0 $     &
 $ X_1 ^{-1} u_{n-1} = u_{n-2 } $  \\
       & $+ \sum_{i=1} ^{ n-2} v_0 ^{-1} a_{n-1-i} u_i +    v_0  u_{n-1}  $ &  for $n\equiv 0\mod 2\mathbb{Z}$;   \\
       & for $n\equiv 0 \mod 2\mathbb{Z} ; $   &  $ X_1 ^{-1} u_{n-1}= \sum_{i=0} ^{n-1} v_0 ^{-1} b_{n-1-i} u_{i} $    \\
       &  $ X_0 ^{-1} \cdot u_{n-1}= u_{n-2} $     &   $+ m_1 u_{n-1} $ \\
       & for $n\equiv 1\mod 2\mathbb{Z}$ ;  &     for $n\equiv 1\mod 2\mathbb{Z}$;   \\

\hline

$E_i \cdot u_j $ & $E_0 \cdot u_0 = \tau_0 u_0  $ ;    &   $E_1 \cdot u_i =0$  \\
                 & $E_0 \cdot u_{2i-1} = \lambda_{2i-1} u_{0}$  ;   &  for any $i$.     \\
                 & $E_0 \cdot u_{2i} = l_0 \lambda_{2i-1} u_{0}  $.   &    \\

\hline

\end{tabular}$\\$

As in \cite{Ma2}, if the set of reflections of a reflection group $W_{\Gamma}$ consists two conjugacy classes then  Marin's generalized Krammer (LK) representation is the direct sum of two components. In above table let $\mathcal{V}_{LK} ^{0} = \Lambda ^{'} <u_0 , u_1 , ..., u_{n-1} >$ be the free $\Lambda ^{'}$ module spanned by $ \{ u_i \} _{i=0,1,...,n-1}$. Operators $X_i ^{\pm} , E_i  \in End ( \mathcal{V}_{LK} ^{0} ),  i=0,1 $ .

 This representation has a $n$ dimensional "mirror" (Definition \ref{defi:anotherevenLK} ) . The sum of these two representations has the right dimension (number of reflections ) , is our version of generalized LK representation of the Artin group $A_{I_2 (2n)}$.

The matrix of $X_0 , X_1$ are as follows when $n$ is even.$\\$

$\left(
  \begin{array}{ccccc}
    l_0  & 0 & 0 & \cdots & 0 \\
    0 & 0 & 1 & \cdots & 0 \\
    m_0 l_0 \lambda_1  & 1 & -m_0  & \cdots & 0 \\
    \vdots  & \vdots & \vdots &  \ddots & \vdots \\
    v_0 ^{-1}a_{n-1} & v_0 ^{-1} a_{n-2} & v_0 ^{-1} a_{n-3} &  \cdots  & v_0 ^{-1} a_0  \\
  \end{array}
\right)_{n\times n},
\left(
  \begin{array}{ccccc}
    0 & 1 & \cdots & 0 & 0 \\
    1 & -m_1  & \cdots & 0 & 0 \\
    \vdots & \vdots & \ddots & \vdots & \vdots \\
    0 & 0 & \cdots & 0 & 1 \\
    0 & 0 & \cdots & 1 & -m_1  \\
  \end{array}
\right)_{n\times n} $

\begin{thm}
\label{thm:evenLK}
The actions in Table 3 defines a representation of the algebra $B_{I_2 (2n)} (\bar{\iota })$ (so also a representation of the Artin group $A_{I_2 (2n)}$). We denote it as $\rho _{LK} ^{0}$.
\end{thm}

\begin{pf} We put the proof  to the last section.
\end{pf} $\\$

The next task is to prove the irreducibility of generalized LK representations and construct invariant bilinear forms.
As in the cases of $I_2 (2n+1)$, we do these by define the following operators.

\begin{defi}
\label{defi:evenBMWspecialprojector}
On the representation space of $\rho_{LK} ^{0}$, we define the following special projectors.
$p_i = [...X_0 X_1 ]_i E_0 [X_1 ^{-1} X_0 ^{-1} ... ]_i $ for $i=0,1,...,n-1$. It is easy to see the image of $p_i$ is the 1 dimensional submodule $\Lambda ^{'} v_i $.
\end{defi}

Next we describe these special projectors explicitly. Suppose
 $$p_{i} (v_j ) = a_{i,j} v_i  .$$
 To describe these coefficients, we only need to consider the cases $i\leq j$, because of the following lemma.

\begin{lem}
\label{lem:easylemmaforprojector}
(1) $E_0 [X_1 ^{-1} X_0 ^{-1} ... ]_{i} [...X_0 X_1 ]_{j} E_0 = a_{i,j} E_0 $;

(2) $a_{j,i} = \overline{a_{i,j}  }$ for any $i,j$.
\end{lem}
\begin{pf}
 For (1),  suppose

 $ E_0 [X_1 ^{-1} X_0 ^{-1} ... ]_{i} [...X_0 X_1 ]_{j} E_0 = a_{i,j} ^{'} E_0 $,

$p_i (v_j ) = [...X_0 X_1 ]_{i} E_0 [X_1 ^{-1} X_0 ^{-1} ... ]_{i} (v_j )
= [...X_0 X_1 ]_{i} E_0 [X_1 ^{-1} X_0 ^{-1} ... ]_{i} [...X_0 X_1 ]_{j} v_0
$

$= a_{i,j} ^{'} [...X_0 X_1 ]_{i} v_{0} = a_{i,j} ^{'} v_i  $.

Compare it with the definition of $a_{i,j}$ we see $a_{i,j} ^{'} = a_{i,j}$.

For (2), to the equality $ E_0 [X_1 ^{-1} X_0 ^{-1} ... ]_{i} [...X_0 X_1 ]_{j} E_0 = a_{i,j} E_0 $ we first apply the automorphism $\Psi$ in Theorem 5.2, which gives $E_0 [X_1  X_0  ... ]_{i} [...X_0 ^{-1} X_1 ^{-1} ]_{j} E_0 =\bar{ a_{i,j} } E_0  $.

Then we apply the anti-automorphism $J$ defined in (5) of remark 5.1, which gives $E_0 [X_1 ^{-1} X_0 ^{-1} ... ]_{j} [...X_0 X_1 ]_{i} E_0 = \overline{ a_{i,j} } E_0   $. So we proved $(2)$.

\end{pf}

\begin{thm}
\label{thm:evenspecialprojectorexplicit}
 (1) $a_{i,i} = \tau_0 $ for any $i$;

(2) $a_{i,j} =l_0 \lambda_{j-i-1}  $ for $i<j$ and $i-j\equiv 0 \mod 2\mathbb{Z}$;

(3) The coefficients $a_{i,j}$ for $i<j$ and $i-j\equiv 1\mod 2\mathbb{Z}$ are determined inductively as follows.

(A) $a_{1, 2i } = \lambda_{2i+1} + m_1 l_0 \lambda_{2i-1} $ for $i\geq 1$;

(B) $a_{2i , 2k+1} = a_{2i-1 , 2k+2 } - m_0 \lambda_{2k+1} \overline{\lambda _{2i-1}}  + m_0 l_0 \lambda_{2k-2i+1} $ for $2i<2k+1 $;

(C) $a_{2i+1 , 2k} = a_{2i, 2k+1} + m_1 l_0 \lambda_{2k-2i-1} $.
\end{thm}

\begin{pf}
First, $(1)$ is true because $E_0 [X_1 ^{-1} X_0 ^{-1} ... ]_{i} [... X_0 X_1 ]_{i} E_0 = E_0 ^{2} = \tau _0 E_0 $.

Then the statement $(2)$ follows from

 $ E_0 [X_1 ^{-1} X_0 ^{-1} ... ]_{i} [... X_0 X_1 ]_{j} E_0
= E_0 [ X_0 ... X_0 X_1 ]_{i-j} E_0 =l_0  E_0 [X_1 ... X_0 X_1 ]_{i-j-1} E_0 = l_0 \lambda_{i-j-1} E_0 .$

Where the first equality sign is because $i<j$ and $i-j\equiv 0 \mod 2\mathbb{Z}$.

For the statement (3), (A) follows from

$E_0 X_1 ^{-1} [X_0 ... X_0 X_1 ]_{2i} E_0 = E_0 (X_1 -m_1 E_1 + m_1  ) [X_0 ... X_0 X_1 ]_{2i} E_0 $

$= E_0 [X_1 ...X_0 X_1 ]_{2i+1} E_0 + m_1 E_0 [X_0 ... X_0 X_1 ]_{2i} E_0
=( \lambda_{2i+1} + m_1 l_0 \lambda_{2i-1} ) E_0 $;

(B) follows from

$E_0 [X_1 ^{-1} X_0 ^{-1} ...X_0 ^{-1} ]_{2i} [X_1 ...X_0 X_1 ]_{2k+1} E_0 =
E_0 [X_1 ^{-1} X_0 ^{-1} ... X_1 ^{-1}  ]_{2i-1} (X_0 -m_0 E_0 +m_0 )[X_1 ...X_0 X_1 ]_{2k+1} E_0
$

$= E_0 [X_1 ^{-1} X_0 ^{-1} ...X_1 ^{-1} ]_{2i-1} [X_0 ...X_0 X_1 ]_{2k+2} E_0 -m_0 \lambda_{2k+1} \bar{ \lambda_{2i-1} } E_0
+ m_0 E_0 [X_0 ...X_0 X_1 ]_{2k-2i+2} E_0
$

$=(a_{2i-1 ,2k+2} -m_0 \lambda_{2k+1} \lambda_{2i-1} + m_0 l_0 \lambda_{2k-2i+1}  ) E_0 $;

(C) follows from

$E_0 [X_1 ^{-1} X_0 ^{-1} ...X_1 ^{-1} ]_{2i+1} [X_0 ...X_1 ]_{2k} E_0 =
E_0 [X_1 ^{-1} X_0 ^{-1} ... X_0 ^{-1} ]_{2i} (X_1 -m_1 E_1 +m_1  ) [X_0 ...X_0 X_1 ]_{2k} E_0
= E_0 [X_1 ^{-1} X_0 ^{-1} ...X_0 ^{-1} ]_{2i} [X_1 ...X_0 X_1 ]_{2k+1} E_0 +
m_1 E_0 [X_0 ... X_0 X_1  ]_{2k-2i} E_0 = ( a_{2i,2k+1} + m_1 l_0 \lambda_{2k-2i}   ) E_0 $.

\end{pf}

\begin{thm}
\label{thm:evenLKirreducible}
As a representation  over the field $F_{\Lambda^{'}}$ , the representation $\rho_{LK} ^{0}$ is irreducible.
\end{thm}

\begin{pf}
Through similar arguments as in the proof of Theorem 3.2, we see if the determinant of the matrix $M= ( a_{i,j} )_{n\times n }$ isn't zero then the representation $\rho_{LK} ^{0}$ is irreducible. So we only need to prove $\det M $ isn't a zero function. For this sake we only need to find a set of values $( l_0 , l_1 , v_0 , v_1  )$ making the determinant nonzero.

Let's have a look at this determinant when:

Fix $l_0 = l_1 =2$; let $v_0 =v_1 =1+\epsilon  $ such that $\epsilon $ is a very small positive number. When $\epsilon $ goes to zero from the positive side, from those relations in theorem 6.2 we easily deduce:

(A) $a_{i,i} = \tau_{0} =\frac{ 2-2 ^{-1}}{m_0 } +1 $ goes to infinity, since $m_1 = m_0 = v_0 - v_0 ^{-1} $ goes to zero from the positive side.

(B) $\lambda_{2i+1}$ goes to $l_0 +1=3$; $ \overline{\lambda_{2i+1}  } $ goes to $l_0 ^{-1} +1 = 1.5 $.

(C) $a_{1,2i}$ goes to $l_0 +1 =3$; more generally  $a_{i,j} $ goes to $l_0 +1 =3 $ if $i<j$ and $i-j\equiv 1 \mod 2\mathbb{Z}$.

(D) Similarly if $i>j$ and $i-j\equiv 1\mod 2\mathbb{Z}$, then  $a_{i,j}$ goes to $l_0 ^{-1} +1=1.5$.

(E) If $i< j$ and $i-j\equiv 0\mod 2\mathbb{Z}$, then $a_{i,j}$ goes to $l_0 (l_0 +1 )= 6$.

(F) If $i>j$ and $i-j \equiv 0\mod 2\mathbb{Z}$, then $a_{i,j}$ goes to $l_0 ^{-1} ( l_0 ^{-1} +1  )= 0.75 $.

So we see when $\epsilon$ goes to $0^{+}$, the diagonal entries of the matrix $M$ goes to infinity, yet all other entries goes to one of numbers in $\{  3,1.5, 6,0.75  \}$. So if $\epsilon $ is small enough, $\det M$ certainly can't be zero. So the statement is proved.

\end{pf}

In the representation $\rho_{LK} ^{0}$ the symmetry in lower indices for "0" and "1" is strongly broken. For example, $E_1$ acts as zero yet $E_1 $ doesn't. In fact applying the symmetry between "0" and "1" in lower indices in Definition 5.1, we can construct another $n$ dimensional representation $\rho_{LK} ^{1}$ as follows.

\begin{defi}
\label{defi:anotherevenLK}
(1) Let the representation space being $\Lambda < \bar{u}_0 ,..., \bar{u}_{n-1}    >$.

(2) As for the action, in Table 4 we replace $( u_i ; X_0 , X_1 , X_0 ^{-1} , X_1 ^{-1} , E_0 , E_1 ; l_0 , l_1 , v_0 , v_1 , m_0 , m_1  )$ with $( \bar{u}_i ; X_1 , X_0 , X_1 ^{-1} , X_0 ^{-1} , E_1 , E_0 ; l_1 , l_0 , v_1 ,v_0 , m_1 ,m_0   )  $ in all places .

\end{defi}

The proof that it indeed produce a representation and  irreducibility of this representation for generic parameters is completely similar, so we omit.

\begin{thm}
\label{thm:evenBMWsemisimpledimension}
As a algebra over the field $F_{\Lambda^{'}}$, the algebra $B_{I_2 (2n)} (\bar{\iota})$ is a $4n+ 2n^2 $ dimensional semisimple algebra.
\end{thm}

\begin{pf}
We have construct two $B_{I_2 (2n)}(\bar{\iota})$ representations: $\rho_{LK} ^{0} $ and $\rho_{LK} ^{1}$. Using the epimorphism $\pi : B_{I_2 (2n)}(\bar{\iota}) \rightarrow H_{I_2 (2n)} (v_0 , v_1 )$ introduced in (3) of Remark 5.1 (where $H_{I_2 (2n)}(v_0 , v_1 )$ is the Hecke algebra ), every (irreducible) representation of $H_{I_2 (2n)}(v_0 , v_1 ) $ provides naturally a (irreducible) $B_{I_2 (2n)}(\bar{\iota} ) $ representation.

Now suppose the data $\bar{\iota} = ( l_0 , l_1 , v_0 , v_1 )$ is generic such that the representations $\rho_{LK} ^{0}$ and $\rho_{LK} ^{1}$ are irreducible and the Hecke algebra $H_{I_2 (2n)}(v_0 , v_1 )$ is semisimple. Suppose the irreducible representations of $H_{I_2 (2n)}(v_0 ,v_1 )$ are $\rho_0 , ..., \rho_{N_{2n } } $. Denote the $B_{I_2 (2n)}(\bar{\iota} )$ representation induced from $\rho_{i}$ as $\bar{\rho}_{i} $. Then since in $\bar{\rho}_i $ ($i=1,..., N_{2n}$) bother $E_0 $ and $ E_1$ act as zero ,and  in $\rho_{LK} ^{0} $ ( $\rho_{LK} ^{1} $   )   only $E_1$ ($ E_0 $) acts as zero, so $\bar{\rho}_1 ,..., \bar{\rho}_{N_{2n}} , \rho_{LK} ^{0} , \rho_{LK} ^{1}$ are different irreducible representations. By Wedderburn-Artin theorem, there is

$\sum_{i=0} ^{ N_{2n}} (\dim \bar{\rho}_{i} )^{2} = \sum_{i=0} ^{N_{2n}} (\dim \rho_{i} )^2 = \dim H_{I_2 (2n)} ( v_0 , v_1) = 4n$.

Using Wedderburn-Artin theorem again we have

$\dim B_{I_2 (2n)}(\bar{\iota }) \geq \sum_{i=0} ^{ N_{2n}} (\dim \bar{\rho}_{i} )^{2} + ( \dim \rho_{LK} ^{0} )^{2}
+ (\dim \rho_{LK} ^{1}  )^{2} = 4n+ 2n^{2} $.

So combine corollary 5.1 we proved for generic parameters $ \dim B_{I_2 (2n)} (\bar{\iota}) = 4n+ 2n^2 $, and at the same time the semisimplicity of $ B_{I_2 (2n)}(\bar{\iota })$.

\end{pf}
$\\$
We can prove the following theorem in the same way as the proof of Theorem 3.3.

\begin{thm}
\label{thm:evenBMWdimensionforanyparameter}
 For any parameters  $\dim B_{I_2 (2n)} ( \bar{\iota}) = 4n+ 2n^{2}$.
\end{thm}

In the same way as the cases $I_2 (2n+1)$ we construct a $A_{I_2 (2n)}$ invariant bilinear form on the representation space $V_{LK} ^{0}$ as follows.  For any element $f(l_0 ,l_1 , v_0 , v_1 ) \in \lambda^{'} $, set $\overline{f} = f(l_0 ^{-1} , l_1 ^{-1} , v_0 ^{-1} , v_1 ^{-1} )$. It is evident the map $f\mapsto \overline{f} $ is a degree 2 automorphism of $\Lambda^{'}$ which may be called conjugation. And we have $\overline{ m_i } = -m_i $.

\begin{defi}
\label{defi:evenLKbilinearform}
Recall the polynomials $a_{i,j}$ defined by those special projectors $p_i$. On the representation space
$V_{LK} ^{0}  = \Lambda^{'} < u_0 , u_1 ,..., u_{n-1}  >$ ,we set
$$ ( u_i , u_j  ) = a_{j,i} = \overline{ a_{i,j} } , $$
then extend it to be a quasi bilinear form $ ( -,-  ): V_{LK} ^{0} \times V_{LK} ^{0} \rightarrow \Lambda^{'} $.
\end{defi}

\begin{thm}
\label{thm:evenbilinearisinvariant}
The bilinear form $(-,- )$ defined in Definition 6.3 is $A_{I_2 (2n)}$ invariant.
\end{thm}
\begin{pf}
We only need to show $( X_k u_i , u_j  ) = ( u_i , X_k ^{-1} u_j ) $ for any $k=0,1$, $0\leq i,j\leq n-1 $. First, because $ p_i ( u_j ) = a_{i,j} u_i  $, so for any $x= \sum_{i=0 } ^{n-1} \mu_j u_j \in B_{G_{2n}}(\bar{\iota}) $,
$ p_i (x) = \sum_{j=0} ^{n-1 } \mu_j a_{i,j} u_i = \sum_{j=0} ^{n-1} \mu_j ( u_j , u_i ) u_i = (x, u_i  ) u_i$. So we have  $ p_j ( X_k u_i ) = (X_k u_i , u_j ) u_j $. From

$ p_j (X_k u_i  )= [...X_0 X_1 ]_{j} E_0 [X_1 ^{-1} X_0 ^{-1} ... ]_{j} X_k [...X_0 X_1 ]_{i} u_0
= (X_k u_i ,  u_j ) [...X_0 X_1 ]_{j} E_0 u_0 $, we see

$E_0 [X_1 ^{-1} X_0 ^{-1} ... ]_{j} X_k [...X_0 X_1 ]_{i} E_0 = ( X_k u_i , u_j   )  E_0 $.

Apply the automorphism $\Psi$ in Theorem 5.2 to both sides get

$E_0 [X_1 X_0 ...  ]_{j} X_k ^{-1} [...X_0 ^{-1} X_1 ^{-1} ]_{i } E_0 = \overline{(X_k u_i , u_j )  } E_0 $.

Then apply the anti-automorphism $J$ in (5) of Remark 5.1 we get

$E_0 [X_1 ^{-1} X_0 ^{-1} ...  ]_{i} X_k ^{-1} [... X_0 ^{-1} X_1 ^{-1}]_{j} E_0 = \overline{ (X_k u_i , u_j  ) } E_0  $.

But at the same time we have
 $E_0 [X_1 ^{-1} X_0 ^{-1} ...  ]_{i} X_k ^{-1} [... X_0 ^{-1} X_1 ^{-1}]_{j} E_0 =(X_k ^{-1} u_j , u_i  ) E_0  $, so we get

  $\overline{ (X_k u_i , u_j )  } = (X_k ^{-1} u_j , u_i )  $, which implies $(X_k u_i , u_j )= (u_i , X_{k} ^{-1} u_j  ) .$

\end{pf}

\section{Cellular structures for $ B_{I_2 (2n)} ( \bar{\iota })$ }

In this section we show the algebra $B_{I_2 (2n)} (\bar{\iota})$ has natural cellular structures. The construction is completely similar with the cases of $B_{I_2 (2n+1)}(m,l) $. Recall there is a natural quotient map
$\pi : B_{I_2 (2n)}(\bar{\iota}) \rightarrow H_{I_2 (2n)} ( v_0 ,v_1 )$ in (3) of Remark 5.1. It is known the Hecke algebra $H_{I_2 (2n)}(v_1 , v_2 )$ has a cellular structure \cite{Ge}. Roughly speaking, the cellular basis of $B_{I_2 (2n)} (\bar{\iota} )$ is a lifting of a cellular basis of $H_{I_2 (2n)}(v_1 , v_2 ) $ together with the set $\{ [...X_0 X_1 ]_{i_1 } E_0 [X_1 X_0 ... ]_{i_2 } |   0\leq i_1 , i_2 \leq n-1 \} \cup  \{ [...X_1 X_0 ]_{i_1} E_1 [X_0 X_1 ... ]_{i_2} | 0\leq i_1 , i_2 \leq n-1 \}$.

Suppose the Hecke algebra $H_{I_2 (2n)}(v_1 ,v_2 )$ has a cellular strucutre $( \Lambda ^{"} , M^{"} , C^{"} , * ) $, as in Definition 4.1. For the anti-automorphism needed for a cellular structure we use the anti-automorphism $J$ introduced in (5) of Definition 5.1. In this section we denote $J(x)$ as $*(x) $.  Denote the ideal of $B_{I_2 (2n)}(\bar{\iota})$ generated by $E_0 , E_1$ as $I_0 $. It is easy to see $I_0 $ is stable under the action of $*$ (or $J$ ), so $J $ induce the anti-automorphism $*$ on the quotient $B_{I_2 (2n)}(\bar{\iota}) / I_0 \cong H_{I_2 (2n)} (v_1 , v_2 )$.

Now the cellular structure $(\Lambda , M , C, *   )  $ of $B_{I_2 (2n)}(\bar{\iota })$ is as follows.

(1) The anti-automorphism $*$ is the $J$ in (5) of Definition 5.1.

(2) $\Lambda = \Lambda^{"} \cup \{ \lambda_{LK} ^{0} , \lambda_{LK} ^{1} \}$. Keep the partial order in $\Lambda^{"}$, and for any $\lambda \in \Lambda^{"}$, let $\lambda_{LK} ^{i}  < \lambda$ for $i=0,1$.

(3) For $\lambda \in \Lambda^{"}$, let $M(\lambda) = M^{"} (\lambda)$ ; and let
$M(\lambda_{LK} ^{i} )=\{  0,1,2,...,n-1 \}  $ for $i=0,1$.

(4) For any $\lambda \in \Lambda^{"}$, choose a total order $"<"$ in $M(\lambda)$. For $S,T\in M(\lambda )$ such that
$ S<T $, set $C_{S,T}$ as any element in $\pi ^{-1} (C^{"} _{S,T} )$.  For $T<S $, set $C_{S,T} = *(C_{T,S} )$. For any $S\in M(\lambda ) $, choose an
element $D_{S} \in \pi^{-1} (C^{"} _{S,S} )$, and set $C_{S,S} = \frac{1}{2} [ D_{S} + *(D_S )  ] $.

(5) For $\lambda_{LK}^{0} $, and $i,j \in M(\lambda_{LK} ^{0})=\{ 0,1,2,...,n-1 \} $, set
$C_{i,j} =[...X_0 X_1 ]_{i} E_0 [X_1 X_0 ... ]_{j}  $;

for $\lambda_{LK} ^{1} $, and $i,j \in M(\lambda_{LK} ^{1})=\{ 0,1,2,..., n-1 \} $, set
$C_{i,j} =[... X_1 X_0  ]_{i}  E_1 [X_0 X_1 ... ]_{j}$.

\begin{thm}
\label{thm:evenBMWcellular}
Above data $(\Lambda , M , C , * ) $ is a cellular structure for $B_{I_2 (2n)}(\bar{\iota }  ) $.
\end{thm}
\begin{pf} The theorem is proved in a similar way with the proof of Theorem \ref{thm:oddBMWcellular}.
\end{pf}

\section{Generalized BMW type algebras for any Coxeter matrixes }

\begin{defi}
\label{defi:generaltypeBMW}
Let $\Gamma=(m_{i,j} )_{n\times n}$ be a Coxeter matrix. The generalized BMW type algebra $B_{\Gamma} (\bar{\iota } )$ is defined as follows.

(1) The data $\bar{\iota} = \{ l_1 , v_1 , l_2 , v_2 ,..., l_n ,v_n \} $, such that, $ l_i = l_j$ and $v_i = v_j $ if
 $ m_{i,j}  $ is odd. Set $m_i = v_i - v_i ^{-1} $ and $\tau_i = \frac{l_i - l_i ^{-1} }{ m_i} +1 $.

 \begin{tabular}{|l|l|l|}
\multicolumn{2}{c}{Table 4: Relations for $B_{\Gamma} (\bar{\iota })$}\\

\hline
   Single $i$  &    (1) $X_i - X_i ^{-1} = m_i (E_i -1 ) $ , (2) $ X_i E_i = E_i X_i = l_i E_i $,     \\
          &   (3) $E_i ^{2} = \tau_i E_i $;   \\

\hline

$m_{i,j} =2$  &  (4) $X_i X_j = X_j X_i $, (5)  $E_i E_j = E_j E_i $, (6) $X_i E_j = E_j X_i  $; \\

\hline

$m_{i,j}= 2l+1$ &   (7) $[X_i X_j ...]_{2l+1}=[X_j X_i ...]_{2l+1} $ ,   (8) $[X_i X_j ... ]_{2l} E_{i} = E_{j} [X_i X_j ...]_{2l} $,  \\
                & (9) $E_i [X_j X_i ... ]_{2k-1}E_i = l_i ^{-1} E_i  $ for $ 1\leq k\leq l $;  \\

\hline

$m_{i,j}=2l $ & (10) $[X_i X_j ... ]_{2l-1} E_j = E_j [X_i X_j ... ]_{2l-1} $, (11) $\Delta_{i,j} E_i = l_i v_i ^{-1} E_i   $,     \\
               & (12) $E_i [X_j X_i ...]_{4k+1} E_i = \frac{m_j}{ m_i}(v_i ^{-1} +l_i ^{-1} )E_i $ for $4k+1<l$, \\
               &(13) $E_i [X_j X_i ...]_{4k+3} E_i = (v_i ^{-1} +l_i ^{-1} )E_i $ for $4k+3<l$, \\
             & (14) $E_i W E_j =0 $ for any word $W$ made from $X_i ^{\pm} , X_j ^{\pm} ,E_i , E_j $. \\

\hline

\end{tabular}$\\$

(2) The set of generators is $\{ X_{i} , X_{i} ^{-1} , E_{i} \} _{i=1,2,...,n } $. For $m_{i,j} $ be even, set $\Delta_{i,j} =[ X_i X_j ... ]_{m_{i,j}}$.

(3) Relations are in the above Table 4.

\end{defi}

For comparison we present the definition of the generalized Brauer type algebra $B_{\Gamma}(\iota)$ defined in \cite{Ch1}.

\begin{defi}
\label{defi:generaltypeBrauer}
Let $\Gamma=(m_{i,j} )_{n\times n}$ be a Coxeter matrix. The type $\Gamma$ Brauer type algebra $Br_{\Gamma}(\iota ) $ is defined as follows.

\begin{tabular}{|l|l|l|}
\multicolumn{2}{c}{Relations for $Br_{\Gamma}(\iota )$ }\\

\hline
     Single $i$ & (1)$S_i ^{2} = 1$, (2) $S_i E_i = E_i S_i = E_i $, (3) $E_i ^{2} =\tau_i E_i$;\\
\hline
     $m_{i,j}=2$ & (4) $S_i S_j =S_j S_i $, (5) $E_i E_j =E_j E_i $, (6) $ S_i E_j = E_j S_i$;\\
\hline
    $m_{i,j}=2l+1$ & (7) $[S_i S_j ... ]_{2l+1} = [S_j S_i ... ]_{2l+1} $,
    (8) $[S_i S_j ...]_{2n} E_i = E_j [S_i S_j ... ]_{2n} $, \\
    &  (9) $E_{i} [S_j S_i ...]_{2k-1} E_{i} = \mu_{ i} E_i $; \\
\hline
    $m_{i,j}= 2l$ & (10) $[S_i S_j ... ]_{2l-1} E_j = E_j [S_i S_j ... ]_{2l-1} = E_j  $, \\
                  & (11) $E_i [S_j S_i ... ]_{4k+1} E_i = (2 \mu_j  ) E_i $ if $l$ is even, and $4k+1<l$, \\
                  & $(11)^{'}$ $E_i [S_j S_i ... ]_{4k+1} E_i = (\mu_i + \mu_j )E_i $ if $l$ is odd $4k+1<l$, \\
                  & (12) $E_i [S_j S_i ... ]_{4k+3} E_i = (2 \mu_i  ) E_i $ if $l$ is even, and $4k+3<l$, \\
                  & $(12)^{'}$ $E_i [S_j S_i ... ]_{4k+3} E_i = (\mu_i + \mu_j )E_i $ if $l$ is odd $4k+3<l$, \\
                  & (13) $E_i W E_j =0$ for any word $W$ mode from $S_i  , E_i , S_j , E_j $.\\

\hline
\end{tabular}

(1) The data $\iota = \{ \tau_i , \mu_i  \}_{i=1,2,...,n }$ such that $\tau_i =\tau_j $ and $\mu_i =\mu_j $ if $m_{i,j} $ is odd.

(2) The set of generators is $\{ S_i , E_i \}_{i=1,2,...,n } $.

(3) The relations is in the following table.
\end{defi}

\section{Proofs}

{\noindent The\ proof\ of\ Theorem \ref{thm:OddLKbasictheorem}}: We need to prove above actions of $X_0 ^{\pm} , X_1 ^{\pm} , E_0 , E_1 $ satisfy all ten relations in Definition \ref{defi:bmwodd} and $X_i X_i ^{-1} = id _{V_{G_{2n+1}}}$.

Denote the operator $[X_0 ... X_1 X_0]_{2n+1}$ and $[X_1 ... X_0 X_1]_{2n+1}$ as $\Delta _1 $ and $\Delta _2$ respectively.

First, the relations $(2) ,(3), (4), (5),(6)$ can be checked by simple computations. To prove the relation (1) is to prove $\Delta _1 = \Delta _2$. By computation of the cases for small $n$'s, we found both $\Delta _1$ and $\Delta _2$  might equal the following operator $\Delta _0 $:
\begin{align*}
 \Delta _0 (v_0 ) &= l v_{2n}  ;\\
\Delta _0 (v_1 )&= ml v_0 + l v_{2n-1} -ml v_{2n};\\
\Delta _0 (v_2 ) &= ml v_1 + l v_{2n-2} -ml v_{2n-1};\\
\Delta _0 (v_3 ) &= ml v_0 + ml v_2 + l v_{2n-3} -ml v_{2n-2} -ml v_{2n };\\
   &......
\end{align*}
\begin{align*}
\Delta _0 (v_{2k+1 }) &= ml \sum _{i=0} ^{k} v_{2i} + l v_{2n-2k-1} -ml \sum_{i=0} ^{k} v_{2n-2i}\quad if\quad  2k+1 \leq n   ;\\
\Delta _0 (v_{2k+1 }) &= ml \sum _{i=0} ^{n-k-1} v_{2i} + l v_{2n-2k-1} -ml \sum_{i=0} ^{n-k-1} v_{2n-2i}\quad if\quad 2k+1 >n £»\\
\Delta _0 (v_{2k }) &= ml \sum _{i=0} ^{k-1} v_{2i+1} + l v_{2n-2k} -ml \sum_{i=0} ^{k-1} v_{2n-2i-1}\quad if\quad  2k\leq n £»\\
\Delta _0 (v_{2k }) &= ml \sum _{i=0} ^{n-k-1} v_{2i+1} + l v_{2n-2k} -ml \sum_{i=0} ^{n-k-1} v_{2n-2i-1}\quad  if \quad 2k>n .
\end{align*}

We present the matrixes of  $2n+1 =5,7$ cases to show what $\Delta_0 $ looks like. $\\$

$\left(
   \begin{array}{ccccc}
    0&0&0&0&l\\
    ml&0&0&l&-ml\\
    0&ml&l&-ml&0\\
    ml&l&0&0&-ml\\
    l&0&0&0&0\\
   \end{array}
 \right)
,\quad  \left(
  \begin{array}{ccccccc}
 0&0&0&0&0&0&l\\
  ml&0&0&0&0&l&-ml\\
 0&ml&0&0&l&-ml&0\\
 ml&0&ml&l&-ml&0&-ml\\
 0&ml&l&0&0&-ml&0\\
 ml&l&0&0&0&0&-ml\\
 l&0&0&0&0&0&0\\
  \end{array}
\right) $ $\\$

Then we prove for any fixed $n$, $\Delta _1 (v_i ) = \Delta _2 (v_i ) =\Delta _0 (v_i )$  by induction on $i$ as follows.

First, when $i=0$,
\begin{align*}
\Delta _1 (v_0 ) &= [X_0 ... X_1 X_0 ]_{2n+1} v_0 =  l v_{2n} = \Delta _0 (v_0 ); \\
\Delta _2 (v_0 ) &=[X_1 ... X_0 X_1 ]_{2n+1} v_0 = X_1 (v_{2n} )= l v_{2n} = \Delta _0 (v_0 ).
\end{align*}

When $i=1$,
\begin{align*}
\Delta _1 (v_1 ) &=[X_0 ... X_0 ]_{2n+1} X_1 v_0 = X_0 X_1 (v_{2n} )= l X_0 (v_{2n} ) \\
&=l(m v_0 - m v_{2n} + v_{2n-1} ) = \Delta _0 (v_1 ) ; \\
\Delta _2 (v_1 ) &= [X_1 ... X_0 X_1  ]_{2n+1} X_1 v_0 = [X_1 ... X_0 ]_{2n} X_1 ^2 v_0 =
[X_1 ... X_0 ]_{2n} (1+ml E_1 - m X_1 ) v_0 \\
&= l [X_1 ... X_1 ]_{2n-1} v_0 + ml v_0 - m X_1 [X_0 ... X_1 ]_{2n} v_0 = l v_{2n-1} +ml v_0 -ml v_{2n} = \Delta _0 (v_1 ) .
\end{align*}
So the statement is true for $i=0,1$. Now assume we have proved the statement for $i<L$. Then we use the following identities to prove the statement for $i=L$ inductively.
\begin{align*}
  (1) \Delta_2 (v_{2k}) &=  [X_1 ... X_1 ]_{2n+1} [X_0 ... X_1 ]_{2k} v_0 =X_1 [X_0 ... X_0 ]_{2n+1} [X_1 ... X_1 ]_{2k-1 } v_0 \\
     &= X_1 [X_0 ... X_1 X_0 ]_{2n+1} v_{2k-1} = X_1 \Delta_1 (v_{2k-1}) .
\end{align*}
\begin{align*}
(2) \Delta_2 (v_{2k+1} ) &=   [X_1 ... X_1 ]_{2n+1} [X_1 ... X_1 ]_{2k+1} v_0 = [X_1 ... X_0 ]_{2n} (1+ml E_1 - m X_1  ) [X_0 ... X_1 ]_{2k} v_0 \\
  &= [X_1 ... X_0 ]_{2n} [X_0 ... X_1 ]_{2k} v_0 + ml v_0 - m[X_1 ... X_1 ]_{2n+1} [X_0 ... X_1 ]_{2k}v_0 \\
   &= X_0 ^{-1} [X_0 ... X_1 X_0 ]_{2n+1} v_{2k} + ml v_0 - m [X_1 ... X_0 X_1 ]_{2n+1} v_{2k}  \\
   &= X_0 ^{-1} \Delta_1 (v_{2k} ) + ml v_0 -m\Delta_2 (v_{2k}) .\\
(3)  \Delta_1 (v_{2k}) & = [X_0 ... X_1 X_0 ]_{2n+1} [X_0 ... X_1 ]_{2k} v_0 \\
    &=  [X_0 ... X_0 X_1 ]_{2n} (1+ml E_0 -m X_0  )[X_1 ... X_1 ]_{2k-1} v_0 \\
     &= X_1 ^{-1} [X_1 ... X_0 X_1]_{2n+1} v_{2k-1} + m v_{2n} - m [X_0 ... X_1 X_0 ]_{2n+1} v_{2k-1} \\
   &= X_1 ^{-1} \Delta_2 (v_{2k-1}) + m v_{2n} -m \Delta_1 (v_{2k-1}) .\\
(4) \Delta_1 (v_{2k+1}) &= [X_0 ... X_1 X_0 ]_{2n+1} [X_1 ... X_0 X_1 ]_{2k+1} v_0 =X_0 [X_1 ... X_0 X_1 ]_{2n+1} [X_0 ... X_0 X_1 ]_{2k} v_0 \\
 &= X_0 [X_1 ... X_0 X_1 ]_{2n+1} v_{2k}  = X_0 \Delta _2 (v_{2k}).
\end{align*}

The induction is divided into the following four cases.

Case $L=2k\leq n+1$. By above equality (1),
\begin{align*}
 \Delta _2 (v_L ) & = X_1 \Delta _1 (v_{2k-1}) =X_1 \Delta _0 (v_{2k-1}) = X_1 ( ml \sum_{i=0} ^{k-1} v_{2i} +l v_{2n-2k+1} -ml \sum_{i=0} ^{k-1} v_{2n-2i} ) \\
 &=ml \sum_{i=0} ^{k-1} v_{2i+1} + l X_1 ^2 v_{2n-2k} -ml \sum_{i=1} ^{k-1} v_{2n-2i+1} -ml^2 v_{2n} \\
&= ml \sum_{i=0} ^{k-1} v_{2i+1} + l v_{2n-2k} - ml \sum_{i=0} ^{k-1} v_{2n-2i-1} = \Delta_0 (v_{2k}) .
\end{align*}

 Where we used

$l X_1 ^2 v_{2n-2k} = l (ml E_1 - m X_1 +1 ) v_{2n-2k} = ml^2 v_{2n} - ml v_{2n-2k +1} + l v_{2n-2k}$.

 By above equality (3) and by induction,
\begin{align*}
 \Delta _1 (v_L ) &= X_1 ^{-1} \Delta _2 v_{2k-1} +m v_{2n} -m \Delta _1 v_{2k-1} = X_1 ^{-1} \Delta _0 v_{2k-1} +m v_{2n} -m \Delta _0 v_{2k-1}\\
& =  X_1 \Delta _0 v_{2k-1} +m v_{2n} -m E_1 \Delta _0 v_{2k-1} =
 X_1 (ml \sum_{i=0} ^{k-1} v_{2i}+ l v_{2n-2k+1} -ml \sum _{i=0} ^{k-1} v_{2n-2i} ) \\
& + mv_{2n} -mE_1
 ( ml \sum_{i=0} ^{k-1} v_{2i} + l v_{2n-2k+1} -ml \sum _{i=0} ^{k-1} v_{2n-2i} ) \\
& = ml \sum_{i=0} ^{k-1} v_{2i+1} + l X_1 ^2 v_{2n-2k} -ml \sum_{i=1} ^{k-1} v_{2n-2i+1} - ml^2 v_{2n}\\
 & + m v_{2n} -m^2 l k v_{2n} +m^2 l (k-1) v_{2n} + m^2 l\tau v_{2n} -ml^2 v_{2n}
 \end{align*}
 \begin{align*}
& = ml \sum_{i=0} ^{k-1} v_{2i+1} +l v_{2n-2k} - ml \sum_{i=0} ^{k-1} v_{2n-2i-1}+ ( - m^2 l +m^2 l \tau +m -ml^2 )  v_{2n}  =  \Delta _0 (v_{2k}).
\end{align*}

 Case $L=2k >n+1$. By using above equality (1),(3) and by induction,
\begin{align*}
\Delta _2 (v_{2k}) & = X_1 \Delta _1 (v_{2k-1}) =X_1 \Delta _0 (v_{2k-1}) = X_1 ( ml \sum_{i=0} ^{n-k} v_{2i} +l v_{2n-2k+1} -ml \sum_{i=0} ^{n-k} v_{2n-2i} ) \\
&= ml \sum_{i=0} ^{n-k} v_{2i+1} + l X_1 ^2 v_{2n-2k} - ml \sum_{i=1} ^{n-k} v_{2n-2i+1} -ml^2 v_{2n} \\
 &=  ml \sum_{i=0} ^{n-k} v_{2i+1} + ml^2 v_{2n} -ml v_{2n-2k+1} +l v_{2n-2k} - ml \sum_{i=1} ^{n-k} v_{2n-2i+1} -ml^2 v_{2n} \\
& = ml \sum_{i=0} ^{n-k-1} v_{2i+1} + lv_{2n-2k} - ml \sum_{i=0} ^{n-k-1} v_{2n-2i-1} = \Delta_0 (v_{2k}) .\\
 \Delta_1 (v_{2k}) &= X_1 \Delta_2 v_{2k-1} + m v_{2n} -m E_1 \Delta _1 v_{2k-1} =
  X_1 (ml \sum_{i=0} ^{n-k} v_{2i} + l v_{2n-2k+1} -ml \sum _{i=0} ^{n-k} v_{2n-2i} ) \\
&+ mv_{2n} -mE_1 ( ml \sum_{i=0} ^{n-k}v_{2i}+ l v_{2n-2k+1} -ml \sum _{i=0} ^{n-k} v_{2n-2i} )\\
& = ml \sum_{i=0} ^{n-k} v_{2i+1} + l X_1 ^2 v_{2n-2k} -ml \sum_{i=1} ^{n-k} v_{2n-2i+1} -ml^2 v_{2n}\\
& + mv_{2n} -m^2 l (n-k+1) v_{2n} -ml^2 v_{2n} + m^2 l (n-k) v_{2n} + m^2 l \tau v_{2n} \\
&= ml \sum_{i=0} ^{n-k-1} v_{2i+1} + lv_{2n-2k} - ml \sum_{i=0} ^{n-k-1} v_{2n-2i-1} + (m- ml^2 -m^2 l +m^2 l \tau ) v_{2n} \\
&= \Delta_0 (v_{2k}) .
\end{align*}

  Case $L=2k+1 \leq n $. By using above equality (2) ,(4) and by induction,
\begin{align*}
\Delta _2 (v_{2k+1} ) &= X_0 [X_0 ... X_1 X_0  ]_{2n+1} v_{2k} +ml v_0 - mE_0 [X_1 ... X_0 X_1 ]_{2n+1} v_{2k} \\
&= X_0 (ml \sum_{i=0} ^{k-1} v_{2i+1} + l v_{2n-2k} -ml \sum_{i=0} ^{k-1} v_{2n-2i-1} )
  + ml v_0 \\
& -m E_0 ( ml \sum_{i=0} ^{k-1} v_{2i+1} + l v_{2n-2k} -ml \sum_{i=0} ^{k-1} v_{2n-2i-1} ) \\
&= ml \sum_{i=0} ^{k-1} v_{2i+2} + l(mv_0 -mv_{2n-2k} + v_{2n-2k-1} ) -ml\sum_{i=0} ^{k-1} v_{2n-2i} + ml v_0 \\
 & -m^2 ll^{-1} kv_0 - mlv_0 + m^2 ll^{-1} k v_0 \\
&= ml\sum_{i=0} ^{k} v_{2i} +lv_{2n-2k-1}-ml\sum_{i=0} ^{k} v_{2n-2i}
  =\Delta_0 (v_{2k+1}).
  \end{align*}
  \begin{align*}
\Delta _1 (v_{2k+1}) &= X_0 [X_1 ...X_0 X_1 ]_{2n+1} v_{2k} = X_0 \Delta_0 v_{2k} =
  X_0 (ml \sum_{i=0} ^{k-1} v_{2i+1} + lv_{2n-2k} -ml \sum_{i=0} ^{k-1} v_{2n-2i-1}  ) \\
& = ml \sum_{i=0} ^{k-1} v_{2i+2} + l (mv_0 -m v_{2n-2k} +v_{2n-2k-1} ) -ml\sum_{i=0} ^{k-1} v_{2n-2i} \\
&= ml\sum_{i=0} ^{k} v_{2i} +lv_{2n-2k-1}-ml\sum_{i=0} ^{k} v_{2n-2i} =\Delta_0 (v_{2k+1}) .
\end{align*}

  Case $L=2k+1 >n $. By using above equality (2) ,(4) and by induction,
 \begin{align*}
\Delta _2 (v_{2k+1} ) &= X_0 [X_0 ... X_1 X_0  ]_{2n+1} v_{2k} +ml v_0 - mE_0 [X_1 ... X_0 X_1 ]_{2n+1} v_{2k} \\
&= X_0 (ml \sum_{i=0} ^{n-k-1} v_{2i+1} + l v_{2n-2k} -ml \sum_{i=0} ^{n-k-1} v_{2n-2i-1} )
  + ml v_0 \\
&-m E_0 ( ml \sum_{i=0} ^{n-k-1} v_{2i+1} + l v_{2n-2k} -ml \sum_{i=0} ^{n-k-1} v_{2n-2i-1} ) \\
&=
  ml\sum_{i=0} ^{n-k-1} v_{2i+2} +l(mv_0 -m v_{2n-2k} +v_{2n-2k+1} ) -ml\sum_{i=0} ^{n-k-1} v_{2n-2i}  \\
&+ mlv_0
   -m^2 ll^{-1} (n-k) v_0 -ml v_0 + m^2 ll^{-1} (n-k) v_0 \\
&=
   ml\sum_{i=0} ^{n-k-1} v_{2i} +lv_{2n-2k-1} -ml\sum_{i=0} ^{n-k-1} v_{2n-2i} = \Delta_0 (v_{2k+1}).  \\
\Delta_1 (v_{2k+1}) &= X_0 [X_1 ...X_0 X_1]_{2n+1} v_{2k} = X_0 \Delta_0 (v_{2k})\\
&=
  X_0 (ml\sum_{i=0} ^{n-k-1} v_{2i+1} +l v_{2n-2k} -ml\sum_{i=0} ^{n-k-1} v_{2n-2i-1}  ) \\
&=
  ml \sum_{i=0} ^{n-k-1} v_{2i+2} +l(mv_0 -mv_{2n-2k} +v_{2n-2k-1} ) - ml\sum_{i=0} ^{n-k-1} v_{2n-2i} \\
& =
   ml\sum_{i=0} ^{n-k-1} v_{2i} +lv_{2n-2k-1} -ml\sum_{i=0} ^{n-k-1} v_{2n-2i} = \Delta_0 (v_{2k+1}).
\end{align*}

So the theorem is proved. $\\$

{\noindent The\ proof\ of\ Theorem \ref{thm:oddspecialprojectorexplicit}}: Since $p_k (v_d ) =E_0 [X_1 ^{-1} X_0 ^{-1} ... ]_k [...X_0 X_1 ]_d E_0 v_k $ by lemma 2.2, to compute these operators, we only need to determine $E_0 [X_1 ^{-1} X_0 ^{-1} ... ]_k [...X_0 X_1 ]_d E_0 $ which we set as $Z_{k,d}$.

Case 1. $( k=d )$   $Z_{k,d} =E_0 ^2 =  \tau E_0 $.

Case 2. $( k-d \equiv 0 \mod 2 , k>d )$  \begin{align*}
Z_{k,d} &= E_0 [X_1 ^{-1} X_0 ^{-1} ... X_0 ^{-1}]_{k-d} E_0 = l^{-1} E_0 [X_1 ^{-1} X_0 ^{-1} ... X_1 ^{-1}]_{k-d-1} E_0 \\
& = l^{-1} l E_0 = E_0 .
\end{align*}

Case 3. $( k-d \equiv 0 \mod 2 , k<d )$
\begin{align*}
Z_{k,d} &= E_0 [X_0 ... X_0 X_1 ]_{d-k} E_0 = l E_0 [X_1 ... X_0 X_1 ]_{d-k-1} E_0 \\
 &= l l^{-1} E_0 = E_0 .
\end{align*}

 Case 4. $( k-d \equiv 1 \mod 2 ,k=0  )$  $ Z_{0,d}= E_0 [X_1 X_0 ...X_1 ]_{d} E_0 = l^{-1} E_0 .$

Case 5. $(k-d\equiv 1\mod 2,  k<2n,  d=2n-1  ) $
\begin{align*}
 Z_{k , 2n-1} &= E_0 [X_1 ^{-1} X_0 ^{-1} ...X_0 ^{-1} ]_{k} [ X_1 ... X_0 X_1 ]_{2n-1 }E_0 =
 E_0 [X_1 ^{-1} ... X_0 ^{-1} X_1 ^{-1} ]_{k-1} X_0 ^{-2} E_1 E_0 \\
 &= E_0 [X_1 ^{-1} ... X_0 ^{-1} X_1 ^{-1} ]_{k-1} (1-ml^{-1} E_0 + m X_0 ^{-1} ) E_1 E_0 =
 l^{-1} E_0 [X_1 ^{-1} X_0 ^{-1} ... X_0 ^{-1} ]_{k-2} E_1 E_0 \\
 &- ml^{-1} E_0 [X_1 ^{-1} ... X_0 ^{-1} X_1 ^{-1}]_{k-1} E_0 + m E_0 [X_1 ^{-1} X_0 ^{-1} ...X_0 ^{-1} ]_k E_1 E_0 \\
 &= l^{-1} E_0 [X_0 ... X_0 X_1 ]_{2n-k+2} E_0 - m l^{-1} l E_0 + m E_0 [X_0 ... X_0 X_1 ]_{2n-k} E_0 = l^{-1} E_0
 \end{align*}

Now we construct a kind of induction process for cases $k-d\equiv 1\mod 2  , k\equiv 0\mod 2, k\geq2$.
\begin{align*}
Z_{k,d} &= E_0 [X_1 ^{-1} X_0 ^{-1}...X_0 ^{-1}]_k [X_0 ^{-1} ... X_1 ^{-1} X_0 ^{-1}]_{2n-d} [X_0 ... X_0 X_1]_{2n} E_0 \\
& = E_0 [X_1 ^{-1} X_0 ^{-1}...X_1 ^{-1}]_{k-1} X_0 ^{-2} [X_1 ^{-1} ... X_1 ^{-1} X_0 ^{-1}]_{2n-d-1} E_1 E_0 \\
&= E_0 [X_1 ^{-1} X_0 ^{-1}...X_1 ^{-1}]_{k-1} ( 1- ml^{-1} E_0 +m X_0 ^{-1} )  [X_1 ^{-1} ... X_1 ^{-1} X_0 ^{-1}]_{2n-d-1} E_1 E_0 \\
&= E_0 [X_1 ^{-1} X_0 ^{-1} ... X_0 ^{-1} ]_{k-2} X_1 ^{-2} [X_0 ^{-1} ... X_1 ^{-1} X_0 ^{-1} ]_{2n-d-2} E_1 E_0 - ml^{-1} l E_0 [X_1 ^{-1} ... X_1 ^{-1} X_0 ^{-1} ]_{2n-d-1} E_0 \\
& + m E_0 [X_1 ^{-1} X_0 ^{-1} ... X_1 ^{-1}]_{k-1} X_0 ^{-1} [X_1 ^{-1} ... X_1 ^{-1} X_0 ^{-1} ]_{2n-d-1} E_1 E_0 \\
&= Term_1 + Term_2 + Term_3 ,
\end{align*}
 where we suppose \begin{align*}
 Term_1 &=  E_0 [X_1 ^{-1} X_0 ^{-1} ... X_0 ^{-1} ]_{k-2} X_1 ^{-2} [X_0 ^{-1} ... X_1 ^{-1} X_0 ^{-1} ]_{2n-d-2} E_1 E_0 ,\\
  Term_2 &=- ml^{-1} l E_0 [X_1 ^{-1} ... X_1 ^{-1} X_0 ^{-1} ]_{2n-d-1} E_0 ,\\
  Term_3 &=  m E_0 [X_1 ^{-1} X_0 ^{-1} ... X_1 ^{-1}]_{k-1} X_0 ^{-1} [X_1 ^{-1} ... X_1 ^{-1} X_0 ^{-1} ]_{2n-d-1} E_1 E_0 .
\end{align*}
Now the second term \begin{align*}
 Term_2 &=  - ml^{-1} l E_0 [X_1 ^{-1} ... X_1 ^{-1} X_0 ^{-1} ]_{2n-d-1}E_0 = -m E_0 [X_0 ... X_0 X_1 ]_{d+1} E_0 \\
  &= -ml E_0 [X_1 ... X_0 X_1 ]_{d} E_0 = -m E_0 .
 \end{align*}

For the third term, if $k<d$, \begin{align*}
  Term_3 &=  m E_0 [X_1 ^{-1} X_0 ^{-1} ... X_1 ^{-1}]_{k-1} X_0 ^{-1} [X_1 ^{-1} ... X_1 ^{-1} X_0 ^{-1} ]_{2n-d-1} E_1 E_0 \\
  &= m E_0 [  X_1 ^{-1} ... X_1 ^{-1} X_0 ^{-1} ]_{2n+k -d-1} E_1 E_0 = m E_0 [X_0 ... X_0 X_1 ]_{d+1-k} E_0 \\
  & = mlE_0 [X_1 ... X_0 X_1 ]_{d-k} E_0 = m E_0 \quad if\quad  k<d  .
 \end{align*}

In case of $k=d+1$, we have

$Term _3 = m E_0 [  X_1 ^{-1} ... X_1 ^{-1} X_0 ^{-1} ]_{2n+k -d-1} E_1 E_0
 = m E_0 [ X_1 ^{-1} ... X_1 ^{-1} X_0 ^{-1}  ]_{2n } [X_0 .... X_1  ]_{2n} E_0 =  m \tau E_0 .$

In case of $k>d+1 $,  we have  \begin{align*}
Term _3 &= m E_0 [  X_1 ^{-1} ... X_1 ^{-1} X_0 ^{-1} ]_{2n+k -d-1} E_1 E_0 = m E_0 [ X_1 ^{-1} ... X_1 ^{-1} X_0 ^{-1}  ]_{2n +k -d-1} [X_0 .... X_1  ]_{2n} E_0 \\
& = m E_0 [X_1 ^{-1} ... X_1 ^{-1} X_0 ^{-1} ]_{k-d-1} E_0 = ml^{-1} E_0 [X_1 ^{-1} ... X_0 ^{-1} X_1 ^{-1} ]_{k-d-2} E_0 = m E_0 .
\end{align*}

 Then we compute the first term.  First if $k\leq d+1 $,
 \begin{align*}
  &Term_1 = E_0 [X_1 ^{-1} X_0 ^{-1} ... X_0 ^{-1} ]_{k-2} X_1 ^{-2} [X_0 ^{-1} ... X_1 ^{-1} X_0 ^{-1} ]_{2n-d-2} E_1 E_0 \\
 &= E_0 [X_1 ^{-1} X_0 ^{-1} ... X_0 ^{-1} ]_{k-2} ( 1-ml^{-1} E_1 +m X_1 ^{-1} ) [X_0 ^{-1} ... X_1 ^{-1} X_0 ^{-1} ]_{2n-d-2} E_1 E_0 \\
 & = E_0 [X_1 ^{-1} X_0 ^{-1} ... X_1 ^{-1} ]_{k-3} X_0 ^{-2} [X_1 ^{-1} ... X_1 ^{-1} X_0 ^{-1} ]_{2n-d-3} E_1 E_0 \\
 &-ml^{-1} E_0 [X_1 ^{-1} X_0 ^{-1} ... X_0 ^{-1} ]_{k-2} E_1 [X_0 ^{-1} ...X_1 ^{-1} X_0 ^{-1} ]_{2n-d-2} E_1 E_0 +m E_0 [X_1 ^{-1} X_0 ^{-1} ... X_1 ^{-1} X_0 ^{-1} ]_{2n+k-d-3}E_1 E_0 \\
 & = E_0 [X_1 ^{-1} X_0 ^{-1} ... X_0 ^{-1}]_{k-2}[X_1 ... X_0 X_1 ]_{d+2} E_0 -ml^{-1} l E_0 [X_0 ... X_0 X_1 ]_{2n-k +2} E_0  + m E_0 [X_0 ... X_0 X_1 ]_{d+3 -k} E_0 \\
&= E_0 [X_1 ^{-1} X_0 ^{-1} ... X_0 ^{-1}]_{k-2}[X_1 ... X_0 X_1 ]_{d+2} E_0 = Z_{k-2 ,d+2 } .
\end{align*}
 If $k= d+3$,
 \begin{align*}
 &Term_1 =E_0 [X_1 ^{-1} X_0 ^{-1} ... X_0 ^{-1}]_{k-2}[X_1 ... X_0 X_1 ]_{d+2} E_0 -ml^{-1} l E_0 [X_0 ... X_0 X_1 ]_{2n-k +2} E_0 \\
 & + m E_0 ^2   =  Z_{k-2 ,d+2} +m (\tau -1) E_0 = Z_{k-2 , d+2 } + (l - l^{-1} ) E_0 .
 \end{align*}

If $k>d+3$,
\begin{align*}
&Term_1 =E_0 [X_1 ^{-1} X_0 ^{-1} ... X_1 ^{-1} ]_{k-3} X_0 ^{-2} [X_1 ^{-1} ... X_1 ^{-1} X_0 ^{-1} ]_{2n-d-3} E_1 E_0 \\
& -ml^{-1} E_0 [X_1 ^{-1} X_0 ^{-1} ... X_0 ^{-1} ]_{k-2} E_1 [X_0 ^{-1} ...X_1 ^{-1} X_0 ^{-1} ]_{2n-d-2} E_1 E_0 +
m E_0 [X_1 ^{-1} X_0 ^{-1} ... X_1 ^{-1} X_0 ^{-1} ]_{2n+k-d-3}E_1 E_0 \\
&=E_0 [X_1 ^{-1} X_0 ^{-1} ... X_0 ^{-1}]_{k-2}[X_1 ... X_0 X_1 ]_{d+2} E_0 -ml^{-1} l E_0 [X_0 ... X_0 X_1 ]_{2n-k +2} E_0  \\
&+ m E_0 [X_1 ^{-1} ... X_1 ^{-1} X_1 ^{-1} ]_{k-d-3} E_0 = Z_{k-2 , d+2} .
\end{align*}

So we have the following inductive equalities.
$$Z_{k,d} = \left\{
  \begin{array}{ll}
    Z_{k-2,d+2}, & k< d; \\
    Z_{k-2, d+2}+(l-l^{-1})E_0 , & k=d+1; \\
    Z_{k-2, d+2} +(l-l^{-1} )E_0, & k=d+3; \\
    Z_{k-2,d+2}, & k>d+3.
  \end{array}
\right.$$
Suppose $k$ is even, and $d$ is odd. By using above computations repeatedly we have

(1) if $k<d$ and $k+d \geq 2n-1$ ,then

$Z_{k,d} = Z_{k-2, d+2} =...=Z_{k+d-2n+1,2n-1} =l^{-1} E_0 $ by case 5.

(2) if $k<d $ and $k+d<2n-1 $ , then

$Z_{k,d} = Z_{k-2,d+2} =...= Z_{0,k+d} =l^{-1} E_0 $ by case 4.

(3) if $k>d+3 $ and $k-d\equiv 3\mod 4$, suppose $k-d=4l+3$, then
\begin{align*}
Z_{k,d} &= Z_{k-2,d+2} = ...=Z_{2l+d+3, 2l+d } = Z_{2l+d+1, 2l+d+2} +(l-l^{-1})E_0 \\
 &= l^{-1} E_0 + (l-l^{-1}  )E_0 = l E_0 \quad by\quad above\quad (1)\quad and\quad (2)\quad ;
\end{align*}

(4) if $k>d+3 $ and $k-d \equiv 1\mod 4 $, suppose $k-d=4l+1$, then
\begin{align*}
Z_{k,d} &= Z_{k-2 ,d+2} = ...= Z_{2l+d+1 , 2l+d  }= Z_{2l+d-1, 2l+d+2}+(l-l^{-1})E_0 \\
&= l^{-1} E_0 + (l-l^{-1} )E_0 = l E_0 \quad by\quad above\quad (1)\quad and\quad (2)\quad ;
\end{align*}

(5) $Z_{d+3 , d} = Z_{d+1, d+2 } +(l-l^{-1})E_0 = l E_0  $;

(6) $Z_{d+1 ,d} = Z_{d-1,d+2} +(l-l^{-1} )E_0 = l E_0 $.

So we have proved the cases $k $ be even and $d$ be odd.  Now the cases $k $ being odd and $d$ being even follows easily from  the following identity
$$ \phi \circ \psi (Z_{k,d}    ) = Z_{d,k }  $$

where the map  $\psi$ and $\phi$ are introduced  in (5) and (3) of  Remark \ref{rem: Gdodd}. So finish the proof. $\\$

{\noindent The\ proof\ of\ Theorem \ref{thm:evenLK}}: We need to show the operators defined in Table 4 satisfy all relations in Definition 5.1. The proof for
relations $4),5),6),8),9),10,11)$ are easy so we omit. The relation $7)$ follows from relation $1) ,2)$. So we only need to check relations $1),2),3),12)$.  First let's see relation $3)$. The cases $X_1 E_1 =l_1 E_1 = E_1 X_1$ and $X_0 E_0 = l_0 E_0$ are evident. So in this case we only need to prove $E_0 X_0 = l_0 E_0$ is satisfied. We prove for any $d$, $E_0 ( X_0 \cdot u_d ) =l_0 E_0 \cdot u_d $.

Case 1 $( d=2i-1<n-1  )$.  $E_0 X_0 u_{2i-1} = E_0 u_{2i} = l_0 \lambda_{2i-1} u_0 $;
$ l_0 E_0 u_{2i-1} = l_0 \lambda_{2i-1} u_0  $.

Case 2 $(d=2i <n-1  )  $.  \begin{align*}
&E_0 X_0 u_{2i} = E_0 ( m_0 l_0 \lambda_{2i-1} u_0 + u_{2i-1} -m_0 u_{2i}   ) \\
&= ( m_0 l_0  \lambda_{2i-1}  \tau_{0} + \lambda_{2i-1} -m_0 l_0 \lambda_{2i-1}   )u_0 =\lambda_{2i-1} ( m_0 l_0 \frac{ l_0 -l_0 ^{-1} }{ m_0 } +1  ) u_0 = l_0 ^{2} \lambda_{2i-1} u_0  ;\\
&l_0 E_0 u_{2i} = l_0 ^{2} \lambda_{2i-1} u_0 .
\end{align*}

Case 3  ( $d= n-1 $ and $n$ is odd ). The check is the same as Case 2.

Case 4  ( $d=n-1 $ and $n$ is even ). Suppose $n=2k$.  $l_0 E_0 u_{n-1} = l_0 \lambda_{2k-1} u_0 $;
\begin{align*}
& E_0 X_0 u_{n-1} = \sum_{i=0} ^{2k-1} v_0 ^{-1} a_{2k-1-i} E_0 u_{i}  =
\sum_{i=0} ^{k-1} v_0 ^{-1} a_{2k-1-2i} E_0 u_{2i} + \sum_{i=1} ^{k} v_0 ^{-1} a_{2k-2i} E_0 u_{2i-1}\\
&= v_0 ^{-1} a_{2k-1} \tau_{0} u_0 + \sum_{i=1} ^{k-1} v_0 ^{-1} l_0 a_{2k-1-2i} \lambda_{2i-1} u_0
+ \sum_{i=1} ^{k} v_0 ^{-1} a_{2k-2i} \lambda_{2i-1} u_0  .
\end{align*}
Coefficient of the  last term  $\sum_{i=1} ^{k} v_0 ^{-1} a_{2k-2i} \lambda_{2i-1} =
 \sum_{i=0} ^{k-1} v_0 ^{-1} a_{2(k-1)-2i} \lambda_{2i+1} = \frac{l_0 ^{-1} + v_0 ^{-1}}{m_0} a_{2k-1}$ by (1) of Lemma 5.2.  For computation of the second term we need to consider the following two cases.

Case 4(1)  ($k=2l+1 $ been odd).   We have \begin{align*}
\sum_{i=1} ^{k-1} l_0 a_{2k-1-2i} \lambda_{2i-1} = \sum_{i=1} ^{2l} a_{4l+1-2i} \lambda_{2i-1} =
\frac{m_1}{m_0}(l_0 + v_0 ) - \frac{l_0 + v_0 }{m_0 } a_{4l+1}
\end{align*}
 by (2) of Lemma 5.2. So  \begin{align*}
 &E_0 X_0 u_{n-1} = v_0 ^{-1} a_{2k-1} \tau_0 u_0 +v_0 ^{-1} [\frac{m_1}{m_0}(l_0 + v_0 )- \frac{l_0 + v_0 }{m_0 }a_{2k-1}  ] u_0 + v_0 ^{-1} \frac{l_0 ^{-1} + v_0 ^{-1} }{m_0} a_{2k-1} u_0 \\
&= l_0 \frac{ m_1}{m_0 } (l_0 ^{-1} + v_0 ^{-1} ) u_0 = l_0 \lambda_{2k-1} u_0 .
\end{align*}
 So this case is proved.

Case 4(2) ($k=2l $ been even  ).  We have \begin{align*}
\sum_{i=1} ^{k-1} l_0 a_{2k-1-2i} \lambda_{2i-1} &= \sum_{i=1} ^{2l-1} l_0 a_{4l-1-2i} \lambda_{2i-1} =
\sum_{i=1} ^{2(l-1) +1 } l_0 a_{4(l-1) +3 -2i} \lambda_{2i-1} \\
&= (l_0 + v_0 ) -\frac{l_0 + v_0  }{ m_0 } a_{2k-1}.
\end{align*} So there is
\begin{align*}
E_0 X_0 u_{n-1}& = v_0 ^{-1} a_{2k-1} \tau_0 + v_0 ^{-1} [(l_0 + v_0  ) - \frac{l_0 +v_0 }{m_0 } a_{2k-1} ]
+ v_0 ^{-1} \frac{l_0 ^{-1} + v_0 ^{-1} }{m_0} a_{2k-1}\\
 &= l_0 (l_0 ^{-1} + v_0 ^{-1} )u_0 = l_0 \lambda_{2k-1} u_0  .
\end{align*}
So finish proof of this case.

The proof for relation $3)$ is also easy so we omit.

At last we prove the key relations $1)$ and $12)$ at the same time.  Set
$$w_i =[... X_0 ^{-1} X_1 ^{-1} ]_{i} u_0 ( = (...(X_0 ^{-1}\cdot (X_1 ^{-1} \cdot u_0 ) )...) , $$
where $X_j ^{-1}$'s are operators defined in  Table 4.  We have proven the operators $X_i , X_i ^{-1} ,E_i  $ $(i=0,1 )$ satisfy relations in $1),7),12)$. Now in the proof of Theorem 5.1 we only use relations other $1), 7), 12)$, which means Theorem 5.1 also holds for those operators in Table 4, equivalently we can prove the following facts in the same way:

$(A) $ $w_{2k} = \sum_{i=0} ^{2k} b_{2k-i} u_i $;

$(B) $ $w_{2k+1} =\sum_{i=0} ^{2k+1} a_{2k+1-i} u_i $.

We need to divide the left proof into two cases, through an  analysis of  iterated action of the operator $X_0 X_1$.

Case (A) ($n=2k$ been even ) From Table 4 we directly see $u_{2i} = (X_0 X_1 )^{i} u_0 $, $ (0\leq i\leq k-1 )$.
Then we have a look at $(X_0 X_1) u_{2k-2}= (X_0 X_1)^{k} u_0 $. First we have $X_1 u_{2k-2} = u_{2k-1} $. Then a check of Table 4 shows $X_0 X_1 u_{2k-1} = X_0 u_{n-1} = \sum_{i=0} ^{2k-1} v_0 ^{-1} a_{2k-1-i} u_i = v_0 ^{-1} w_{2k-1} $ by above  $(B) $. So we have proved $[...X_0 X_1 ]_{n} u_0 = v_0 ^{-1} [...X_0 ^{-1} X_1 ^{-1} ]_{n-1} u_0   $,
which implies $[...X_0 X_1]_{2n-1} u_0 = v_0 ^{-1} u_0 $, which implies further $\Delta u_0 =[...X_1 X_0 ]_{2n} u_0 = l_0 v_0 ^{-1} u_0 $.

So we have $(X_0 X_1) u_{n-2} = v_0 ^{-1} w_{n-1 }$. Repeat the action of $X_0 X_1$ again we have

$(X_0 X_1)^{k+1} u_0 = (X_0 X_1 )^{2} u_{n-2}= v_0 ^{-1} w_{n-3} ,..., (X_0 X_1) ^{n-1} u_0 =v_0 ^{-1} w_{1}
= v_0 ^{-1} X_1 ^{-1} u_0 $, and $(X_0 X_1 )^{n} u_0 = l_0 v_0 ^{-1} u_0 $.  Which form a " cycle" but the result is a multiplication of the constant $l_0 v_0 ^{-1}$. The elements appear in the cycle is

$u_0 , u_2, u_4 ,...,u_{n-2} (= u_{2k-2}) , w_{n-1} , w_{n-3} ,..., w_1 .$\ \ \   $( * )$

The action of $X_0 X_1$ send every element in this sequence to a constant multiple of the next one ($X_0 X_1 w_1 = l_0 u_0 $  ) . So for any element $X$ in this sequence we have the same identity $\Delta X= l_0 v_0 ^{-1} X $.

Next we see the transformation matrix from the basis $u_0 , u_1 , ..., u_{n-1}  $ to

 $u_0 , w_1 , u_2 , w_3 ,..., u_{n-2} , w_{n-1} $ is a upper triangular matrix whose diagonal entries are all 1. So the sequence $(\alpha )$ is also a basis of the representation space, thus we have $\Delta = l_0 v_0 ^{-1} I_{n} $. So we have

$[...X_1 X_0 ]_{2n} = X_1 [... X_0 X_1 ]_{2n} X_1 ^{-1} = l_0 v_0 ^{-1} I_{n} = [...X_0 X_1 ]_{2n} $.

Relation $12)$ follows immediately also.

Case (B) ($n=2k+1$ been odd ) Proof this this case is similar so we omit the details. We still prove it by analysis the action of $X_0 X_1 $. This time it produce a sequence

$u_0 , u_2 ,..., u_{n-1} , w_{n-2} ,w_{n-4} ,... w_1  $.  \ \ \   $( ** )$

Where we have $(X_0 X_1 )^{i} u_0 = u_{2i} $ $(0\leq i\leq k)$, $(X_0 X_1  )^{k+i} u_0 = v_0 ^{-1} w_{n-2i } $ $(1\leq i\leq k )$, $\Delta u_0 = (X_0 X_1 )^{n} u_0 = l_0 v_0 ^{-1} u_0 $.

For the same reason the  elements in $(**)$ still form a basis $u_0 , w_1 , u_2 , w_3 ,...,w_{n-2} , u_{n-1} $. So we have proved relation $1) ,12)$ as well in this case.

\hspace{-0.70cm} {\sc School of Mathematics}

\nd{\sc Hefei university of technology}

\nd {\sc Hefei 230009 China}

\nd {\sc E-mail addresses}: {\sc Zhi Chen} ({\tt
zzzchen@ustc.edu.cn}).

\end{document}